\documentclass[12pt,oneside]{amsart}

\usepackage[american]{babel}
\usepackage[%
  style=authoryear-comp,
  backend=bibtex,
  eprint=false,
  doi=false,
  url=false,
  firstinits=true,
  bibstyle=jltbibstyle]{biblatex}
\bibliography{patents,taffy}

\usepackage{ifthen}
\def\arxivversion{true}
\newcommand{\ifarxiv}[2]{\ifthenelse{\equal{\arxivversion}{true}}{#1}{#2}}
\def\hiresversion{true}
\newcommand{\ifhires}[2]{\ifthenelse{\equal{\hiresversion}{true}}{#1}{#2}}

\usepackage{amsmath,amsfonts,amssymb,amsthm,txfonts,pxfonts,amscd}
\usepackage[mathcal]{euscript}
\usepackage{graphicx}
\usepackage{subfigure}
\usepackage[letterpaper]{geometry}
\usepackage[T1]{fontenc}
\usepackage{hyperref}
\usepackage{tikz}
\usetikzlibrary{decorations.pathreplacing,decorations.markings,shapes.misc}
\usetikzlibrary{calc}
\usetikzlibrary{fadings}

\definecolor{lgray}{rgb}{.85,.85,.85}

\tikzset{
  on each segment/.style={
    decorate,
    decoration={
      show path construction,
      moveto code={},
      lineto code={
        \path [#1]
        (\tikzinputsegmentfirst) -- (\tikzinputsegmentlast);
      },
      curveto code={
        \path [#1] (\tikzinputsegmentfirst)
        .. controls
        (\tikzinputsegmentsupporta) and (\tikzinputsegmentsupportb)
        ..
        (\tikzinputsegmentlast);
      },
      closepath code={
        \path [#1]
        (\tikzinputsegmentfirst) -- (\tikzinputsegmentlast);
      },
    },
  },
  mid arrow/.style={postaction={decorate,decoration={
        markings,
        mark=at position .6 with {\arrow[#1]{>}}
      }}},
  mid arrow2/.style={postaction={decorate,decoration={
        markings,
        mark=at position .6 with {\arrow[#1]{>>}}
      }}},
  mid arrow3/.style={postaction={decorate,decoration={
        markings,
        mark=at position .6 with {\arrow[#1]{latex}}
      }}}
}

\tikzset{->-/.style={decoration={
  markings,
  mark=at position #1 with {\arrow{>}}},postaction={decorate}}}
\tikzset{->>-/.style={decoration={
  markings,
  mark=at position #1 with {\arrow{>>}}},postaction={decorate}}}
\tikzset{-|>-/.style={decoration={
  markings,
  mark=at position #1 with {\arrow{latex}}},postaction={decorate}}}

\newcommand{\tikzbgellipsefour}{%
  \fill[lgray] ($(p1) + .5*(\tw,0)$) ellipse (2.5cm and 1.5cm);
}

\newcommand{\tikzbgellipsefive}{%
  \fill[lgray] (q2) ellipse (3.5cm and 1.5cm);
}

\newcommand{\tikzbgellipsesix}{%
  \fill[lgray] (O) ellipse (4.75cm and 1.5cm);
}

\usepackage{units}

\ifarxiv{
  \ifhires{
    \graphicspath{{figs/}{figs/hi/}{figs_extra_arxiv/}{figs_extra_arxiv/hi/}}
  }{
    \graphicspath{{figs/}{figs/lo/}{figs_extra_arxiv/}{figs_extra_arxiv/lo/}}
  }
}{
  \ifhires{
    \graphicspath{{figs/}{figs/hi/}}
  }{
    \graphicspath{{figs/}{figs/lo/}}
  }
}

\newcommand{\keepnote}[1]{}

\newcommand{\mathnotation}[2]{\newcommand{#1}{\ensuremath{#2}}}

\newcommand{\hvect}[2]{(\,#1 \ \, #2\,)}
\newcommand{\vect}[2]{\hvect{#1}{#2}}
\newlength{\subfigwidth}
\renewcommand{\l}{\left}                        % \left
\renewcommand{\r}{\right}                       % \right
\mathnotation{\id}{\mathrm{id}}              % identity map
\mathnotation{\torus}{T^2}                   % torus
\renewcommand{\S}{S}                         % a surface
\mathnotation{\Sq}{S}                        % torus mod hyperelliptic inv
\mathnotation{\pW}{p}                        % fixed pts of hyperelliptic inv
\mathnotation{\SL}{\mathrm{SL}}              % Special Linear group
\mathnotation{\Z}{\mathbb{Z}}                % integers

\begin{document}

\title{The mathematics of taffy pullers}
\author{Jean-Luc Thiffeault}
\thanks{Supported by NSF grant CMMI-1233935%
  \ifarxiv{.  Contains an extra appendix compared to the version
    published in
    \href{http://dx.doi.org/10.1007/s00283-018-9788-4}{\emph{Mathematical
    Intelligencer}}}{}}
\address{Department of Mathematics,
  University of Wisconsin, Madison, WI 53706, USA}
\email{jeanluc@math.wisc.edu}

\begin{abstract}
  We describe a number of devices for pulling candy, called taffy pullers,
  that are related to pseudo-Anosov maps of punctured spheres. Though the
  mathematical connection has long been known for the two most common taffy
  puller models, we unearth a rich variety of early designs from the patent
  literature, and introduce a new one.
\end{abstract}

\maketitle

\section*{Introduction}

\begin{figure}
\begin{center}
  \includegraphics[width=.7\textwidth]{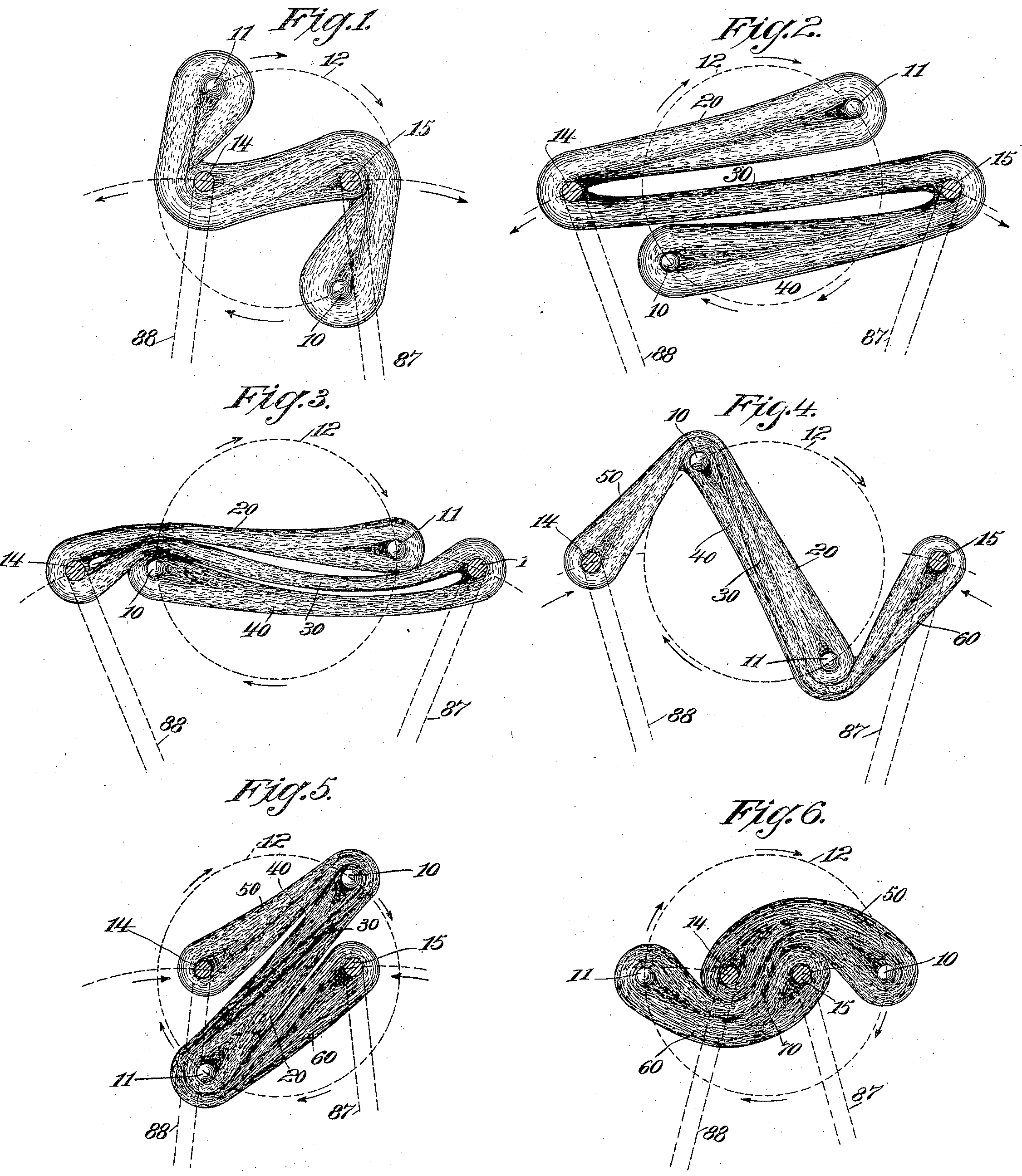}
\end{center}
\caption{Action of the taffy puller patented by
  \textcite{richards_process_1905}.  The rod motion is equivalent (conjugate)
  to that of Fig.~\ref{fig:Thibodeau1903}.}
\label{fig:Richards1905}
\end{figure}

Taffy is a type of candy made by first heating sugar to a critical
temperature, letting the mixture cool on a slab, then repeatedly `pulling' ---
stretching and folding --- the resulting mass.  The purpose of pulling is to
get air bubbles into the taffy, which gives it a nicer texture.  Many devices
have been built to assist pulling, and they all consist of a collection of
fixed and moving rods, or pins.  Figure~\ref{fig:Richards1905} shows the
action of such a taffy puller from an old patent.  Observe that the taffy
(pictured as a dark mass) is stretched and folded on itself repeatedly.  As
the rods move, the taffy is caught on the rods and its length is forced to
grow exponentially.  The effectiveness of a taffy puller is directly
proportional to this growth, since more growth implies a more rapid trapping
of the air bubbles.  Given a pattern of periodic rod motion, regarded as
orbits of points in the plane, the mathematical challenge is to compute the
growth.

We will describe in broad terms how the computation of growth is achieved.
The framework involves the topological dynamics ideas pioneered by William
Thurston, but we will shy away from a complete treatment involving rigorous
definitions.  Instead we will boil down the computation to its essence: the
relationship between maps of the torus and those on a punctured sphere.
Computations on the former involve simple linear algebra, and the taffy
pullers are described by the latter.  We will also show an explicit example
that involves surfaces of higher genus than a torus, which allow us to
describe taffy pullers with more than three or four moving rods.  Throughout,
we will give examples of taffy pullers from the patent literature as well as a
newly-invented one.  Finally, we answer the question: which taffy puller is
the `best' in a mathematical sense?

\section*{Some history}

Until the late 19th century, taffy was pulled by hand --- an arduous task.
The process was ripe for mechanization.  The first patent for a mechanical
taffy puller was by \textcite{firchau_machine_1893}: his design consisted of
two counter-rotating rods on concentric circles.  This is not a `true' taffy
puller: a piece of taffy wrapped around the rods will not grow exponentially.
Firchau's device would have been terrible at pulling taffy, but it was likely
never built.

In 1900, Herbert M.\ Dickinson invented the first nontrivial taffy puller, and
described it in the trade journal \emph{The Confectioner}.  His machine
involved a fixed rod and two rods that move back-and-forth.  The moving rods
are `tripped' to exchange position when they reach the limit of their motion.
Dickinson later patented the machine \autocite{dickinson_candy-pulling_1906}
and assigned it to Herbert L.\ Hildreth, the owner of the Hotel Velvet on Old
Orchard Beach, Maine.  Taffy was especially popular at beach resorts, in the
form of salt water taffy (which is not really made using salt water).
Hildreth sold his `Hildreth's Original and Only Velvet Candy' to the Maine
tourists as well as wholesale, so he needed to make large quantities of taffy.
Though he was usually not the inventor, he was the assignee on several taffy
puller patents in the early 1900s.  In fact several such patents were filed in
a span of a few years by several inventors, which led to lengthy legal
wranglings.  Some of these legal issues were resolved by Hildreth buying out
the other inventors; for instance, he acquired one patent for \$75,000 (about
two million of today's dollars).  Taffy was becoming big business.

Shockingly, the taffy patent wars went all the way to the US Supreme Court.
The opinion of the Court was delivered by Chief Justice William Howard Taft.
The opinion shows a keen grasp of topological dynamics (\emph{Hildreth v.\
  Mastoras}, 1921):
\begin{quote}
  The machine shown in the Firchau patent
  % comprises two discs which are
  % rotated in opposite directions. On each disc is a finger which projects
  % into
  % a drum, into which the candy is put. The
  [has two pins that] pass each other twice during each revolution
  % of the disc
  [\dots] and move in concentric circles, but do not have the relative
  in-and-out motion or Figure 8 movement of the Dickinson machine. With only
  two hooks there could be no lapping of the candy, because there was no third
  pin to re-engage the candy while it was held between the other two pins. The
  movement of the two pins in concentric circles might stretch it somewhat and
  stir it, but it would not pull it in the sense of the art.
  % The Firchau
  % device never, so far as appears in the record, made candy experimentally
  % or
  % otherwise. Indeed, no candy was commercially pulled by machine before or
  % after the issuing of the Firchau patent in 1893 until the introduction of
  % the Dickinson principle, nine or ten years later.
\end{quote}
The Supreme Court opinion displays the fundamental insight that at least three
rods are required to produce some sort of rapid growth.  Moreover, the `Figure
8' motion is identified as key to this growth.  We shall have more to say on
this rod motion as we examine in turn the different design principles.

\begin{figure}
\begin{center}
\subfigure[]{
  \includegraphics[width=.5\textwidth]{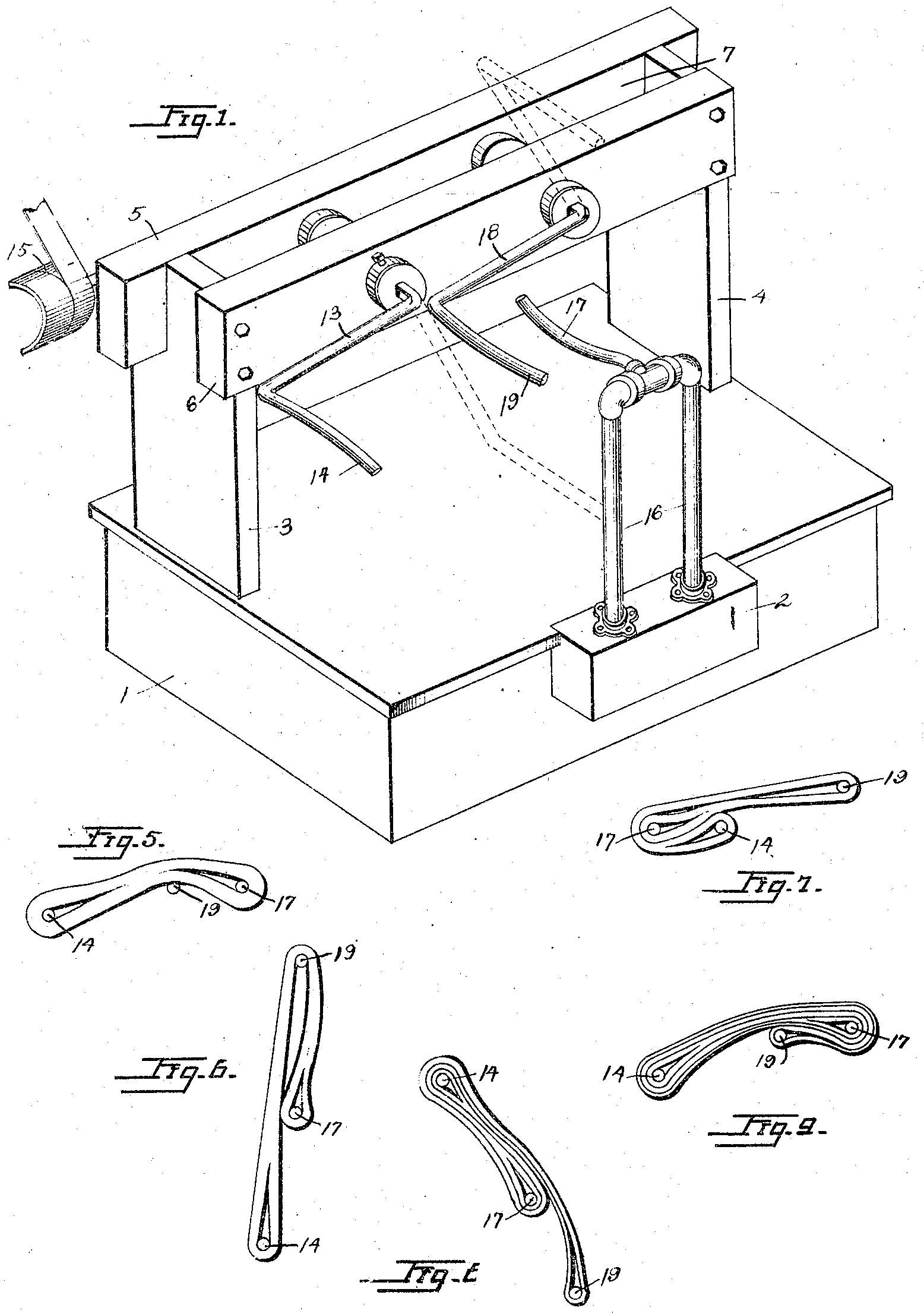}
  \label{fig:Robinson1908_device}
}\hspace{1em}
\subfigure[]{
  \raisebox{1.35in}{
  \includegraphics[height=.2\textheight]{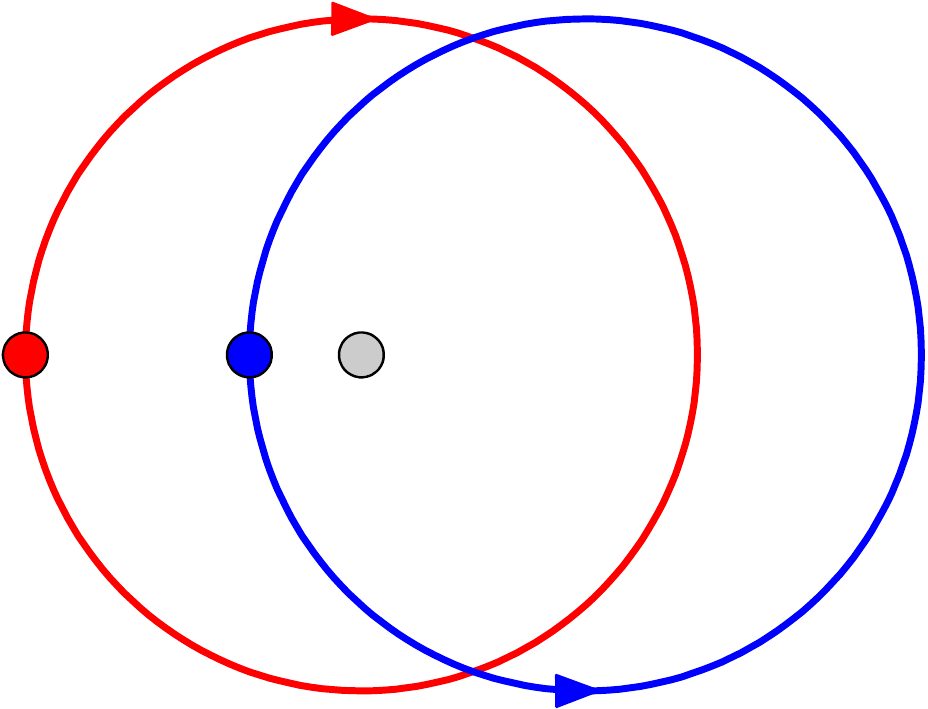}}
  \label{fig:Robinson1908_rodmotion}
}
\end{center}
\caption{(a) Taffy puller from the patent of
  \textcite{robinson_candy-pulling_1908}.  (b) The motion of the rods.}
\label{fig:Robinson1908}
\end{figure}

The Dickinson taffy puller may have been the first, but it was overly
complicated and likely never used to make large quantities of candy.  A
similar rod motion can be obtained by a much simpler mechanism, which was
introduced in a patent by \textcite{robinson_candy-pulling_1908} and is still
in use today.  In this device, two rods move in counter-rotating orbits around
a fixed rod (Fig.~\ref{fig:Robinson1908}).  We call this design the
\emph{standard 3-rod taffy puller}.

\section*{Three-rod taffy pullers}

\begin{figure}
\begin{center}
  \begin{minipage}{.3\textwidth}
    \subfigure[]{
      \includegraphics[width=\textwidth]{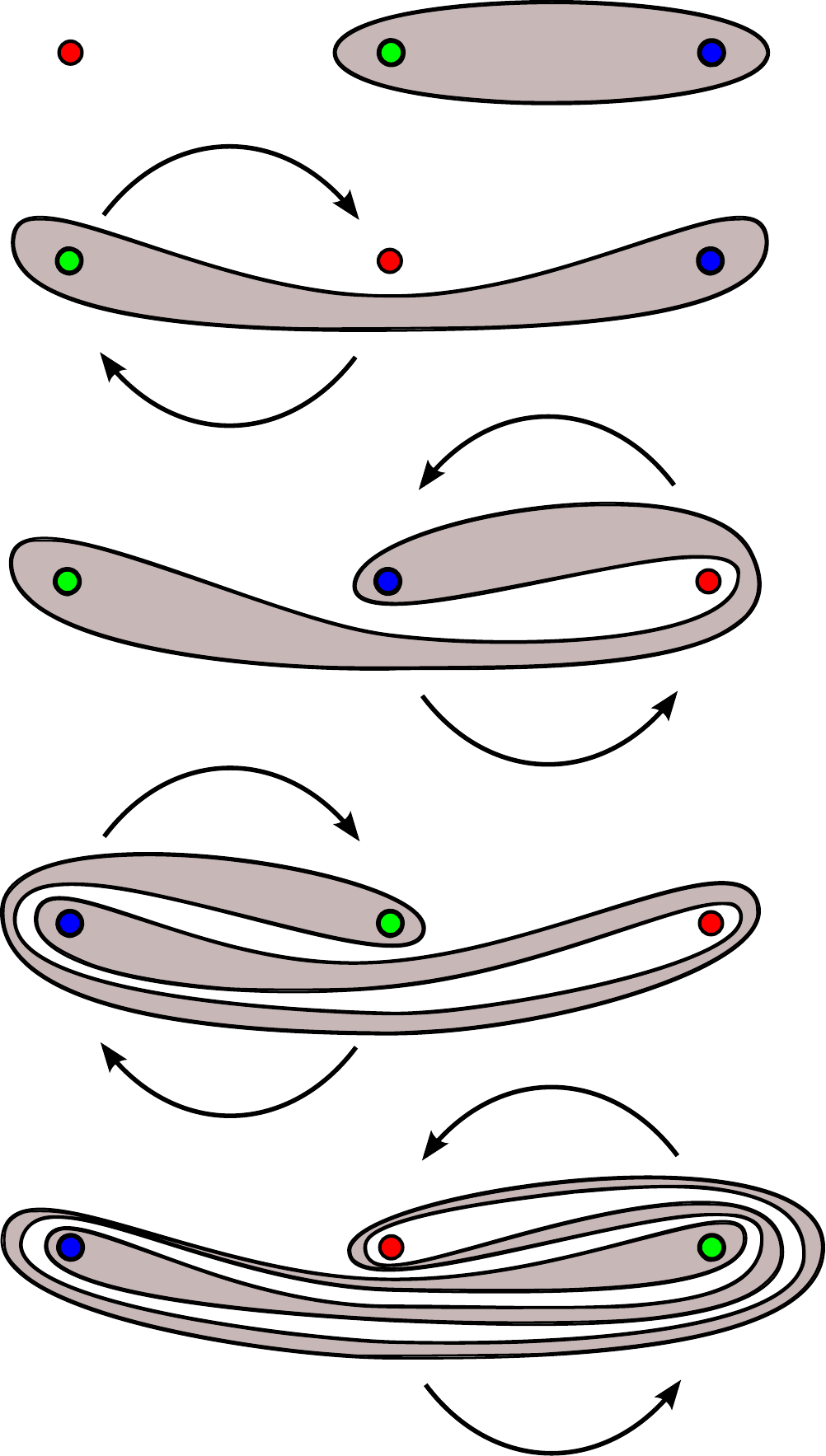}
      \label{fig:s1s-2_loop}
    }
  \end{minipage}
  \hspace{2em}
  \begin{minipage}{.4\textwidth}
    \subfigure[]{
      \includegraphics[width=\textwidth]{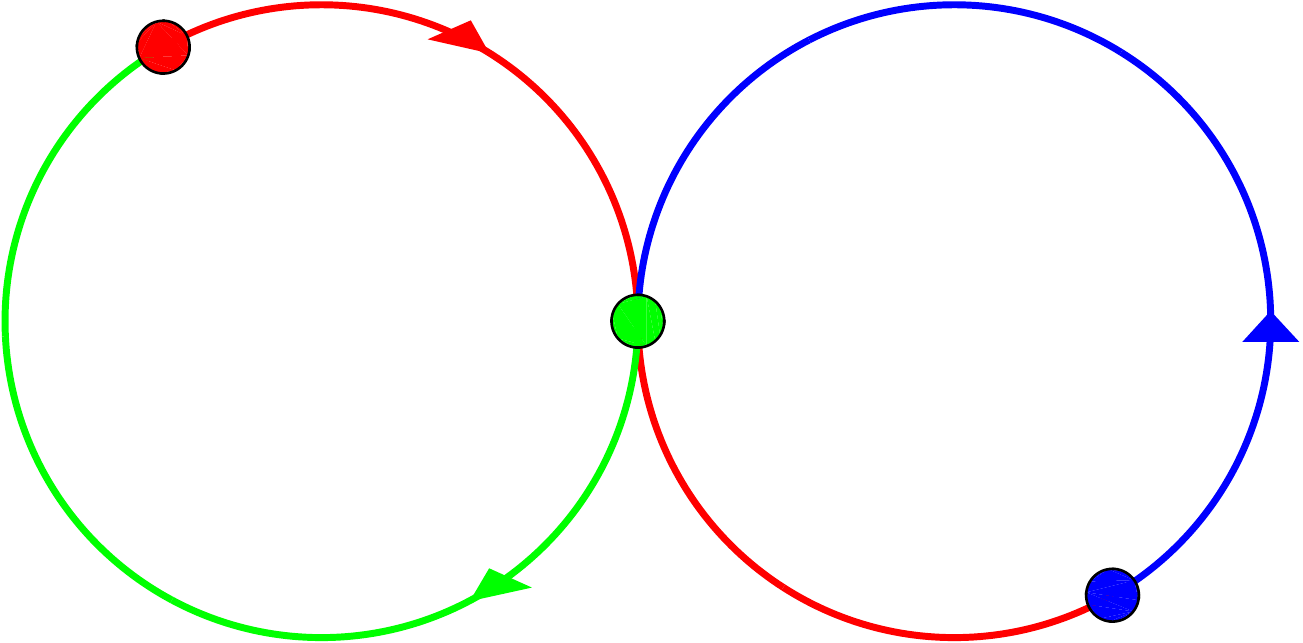}
      \label{fig:Nitz1918_rodmotion}
    }

    \subfigure[]{
      \includegraphics[width=\textwidth]{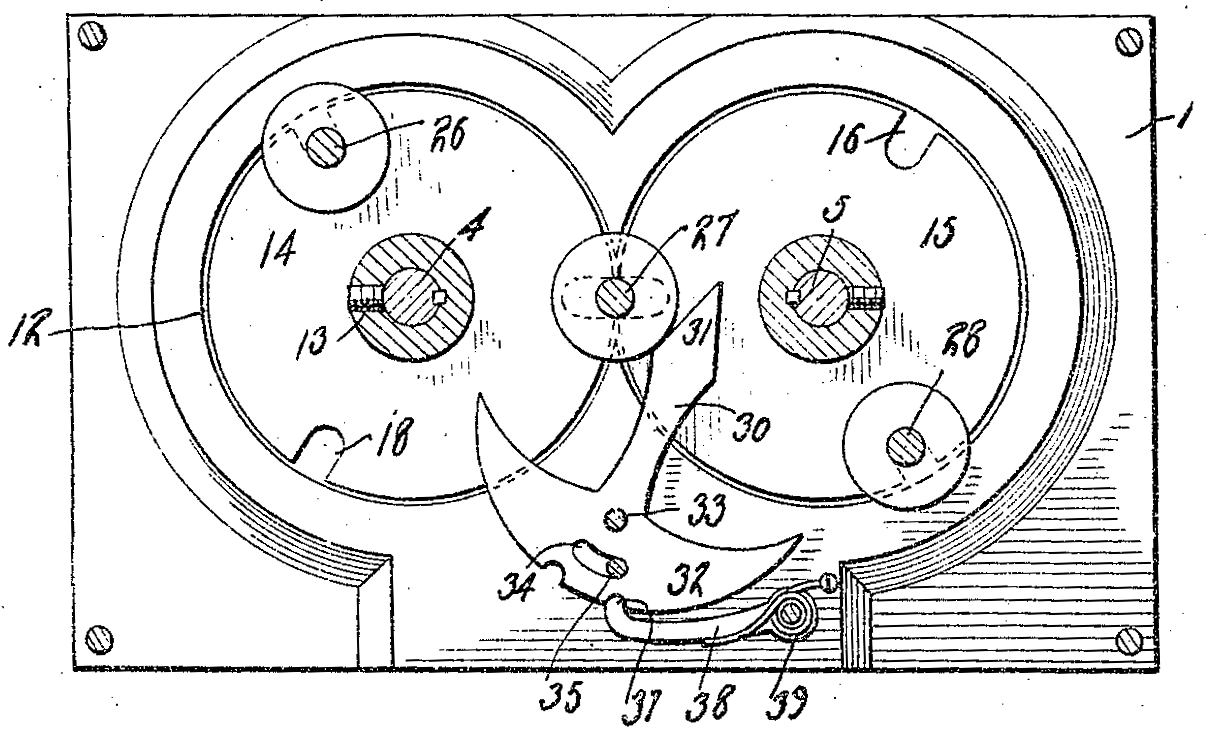}
      \label{fig:Nitz1918_device}
    }
  \end{minipage}
\end{center}
\caption{(a) The action of a 3-rod taffy puller.  The first and second rods
  are interchanged clockwise, then the second and third rods are interchanged
  counterclockwise.  (b) Each of the three rods moves in a Figure-8.  (c)
  Taffy puller from the patent of \textcite{nitz_candy-puller._1918}, where
  rods alternate between the two wheels.}
\label{fig:s1s-2}
\end{figure}

Taffy pullers involving three rods (some of which may be fixed) are the
easiest to describe mathematically.  The action of arguably the simplest such
puller, from the mathematical standpoint, is depicted in
Fig.~\ref{fig:s1s-2_loop}.  By action, we mean the effect of the puller on a
piece of abstract `taffy.'  For this puller, the first and second rods are
interchanged clockwise, then the second and third are interchanged
counterclockwise.
%Here first, second, and third refer to the current
%\emph{position} of a rod from left to right, not to the rod itself.
Notice that each rod undergoes a `Figure 8' motion, as shown in
Fig.~\ref{fig:Nitz1918_rodmotion} for~$\tfrac13$ period.

\begin{figure}
\begin{center}
  \subfigure[]{
    \resizebox{.44\subfigwidth}{!}{\begin{tikzpicture}
  \newcommand{\tw}{36pt}
  \begin{scope}[line width=.7pt]
    \draw[postaction={on each segment={mid arrow}}]
      (0,0) -- (0,1/2*\tw);
    \draw[postaction={on each segment={mid arrow}}]
      (\tw,0) -- (\tw,1/2*\tw);
    \draw[postaction={on each segment={mid arrow}}]
      (0,\tw) --(0,1/2*\tw);
    \draw[postaction={on each segment={mid arrow}}]
      (\tw,\tw) -- (\tw,1/2*\tw);
    \draw[postaction={on each segment={mid arrow2}}]
      (0,0) -- (1/2*\tw,0);
    \draw[postaction={on each segment={mid arrow2}}]
      (0,\tw) -- (1/2*\tw,\tw);
    \draw[postaction={on each segment={mid arrow2}}]
      (\tw,0) -- (1/2*\tw,0);
    \draw[postaction={on each segment={mid arrow2}}]
      (\tw,\tw) -- (1/2*\tw,\tw);
  \end{scope}
  \begin{scope}[line width=.5pt]
    \draw[red,postaction={on each segment={mid arrow}}]
      (0,1/2*\tw) -- (1/2*\tw,1/2*\tw);
    \draw[red,postaction={on each segment={mid arrow}}]
      (\tw,1/2*\tw) -- (1/2*\tw,1/2*\tw);
    \draw[red,postaction={on each segment={mid arrow3}}]
      (1/2*\tw,\tw) -- (1/2*\tw,1/2*\tw);
    \draw[red,postaction={on each segment={mid arrow3}}]
      (1/2*\tw,0) -- (1/2*\tw,1/2*\tw);
  \end{scope}
  \fill[gray,nearly transparent]
    (0,0) -- (1/2*\tw,0) -- (1/2*\tw,1/2*\tw) -- (0,1/2*\tw) -- cycle;
  \fill[gray,nearly transparent, shift={(1/2*\tw,1/2*\tw)}]
    (0,0) -- (1/2*\tw,0) -- (1/2*\tw,1/2*\tw) -- (0,1/2*\tw) -- cycle;
  \node at (1/4*\tw,1/4*\tw) {$A$};
  \node[rotate=180] at (3/4*\tw,3/4*\tw) {$A$};
  \node at (3/4*\tw,1/4*\tw) {$B$};
  \node[rotate=180] at (1/4*\tw,3/4*\tw) {$B$};
\end{tikzpicture}}
    \label{fig:torus-hyper}
  }
  \hspace{2em}
  \subfigure[]{
    \resizebox{.5\subfigwidth}{!}{\begin{tikzpicture}
  \newcommand{\tw}{36pt}
  \begin{scope}[line width=.7pt]
    \draw[postaction={on each segment={mid arrow}}]
      (0,0) -- (0,1/2*\tw);
    \draw[postaction={on each segment={mid arrow}}]
      (\tw,0) -- (\tw,1/2*\tw);
    \draw[postaction={on each segment={mid arrow2}}]
      (0,0) -- (1/2*\tw,0);
    \draw[postaction={on each segment={mid arrow2}}]
      (\tw,0) -- (1/2*\tw,0);
  \end{scope}
  \begin{scope}[line width=.5pt]
    \draw[red,postaction={on each segment={mid arrow}}]
      (0,1/2*\tw) -- (1/2*\tw,1/2*\tw);
    \draw[red,postaction={on each segment={mid arrow}}]
      (\tw,1/2*\tw) -- (1/2*\tw,1/2*\tw);
    \draw[red,postaction={on each segment={mid arrow3}}]
      (1/2*\tw,0) -- (1/2*\tw,1/2*\tw);
  \end{scope}
  \fill[gray,nearly transparent]
    (0,0) -- (1/2*\tw,0) -- (1/2*\tw,1/2*\tw) -- (0,1/2*\tw) -- cycle;
  \node at (1/4*\tw,1/4*\tw) {$A$};
  \node at (3/4*\tw,1/4*\tw) {$B$};
  \filldraw[black] (0,0) circle (.02*\tw)
    node[below,scale=.75] {\tiny $\pW_0$};
  \filldraw[black] (1/2*\tw,0) circle (.02*\tw)
    node[below,scale=.75] {\tiny $\pW_1$};
  \filldraw[black] (0,1/2*\tw) circle (.02*\tw)
    node[above,scale=.75] {\tiny $\pW_3$};
  \filldraw[black] (1/2*\tw,1/2*\tw) circle (.02*\tw)
    node[above,scale=.75] {\tiny $\pW_2$};
\end{tikzpicture}}
    \label{fig:torus-hyper2}
  }
  \hspace{1em}
  \subfigure[]{
    \raisebox{.5em}{
    \resizebox{!}{.6\subfigwidth}{\begin{tikzpicture}
  \newcommand{\tw}{36pt}
  \coordinate (p0) at (-.7*\tw,0*\tw);
  \coordinate (p1) at (-.2*\tw,-.15*\tw);
  \coordinate (p2) at (.3*\tw,-.15*\tw);
  \coordinate (p3) at (.8*\tw,0*\tw);
  \coordinate (p01c) at ($(p0)!.4!(p1) + (0,-.05)$);
  \coordinate (p12c) at ($(p1)!.5!(p2) + (0,-.025)$);
  \coordinate (p32c) at ($(p3)!.5!(p2) + (0,-.05)$);
  \begin{scope}[line width=.7pt]
  \draw[] (0,0) circle (\tw);
  \shadedraw[shading=ball,ball color=white]
    (0,0) circle (\tw);
  \end{scope}
  \begin{scope}[line width=.5pt]
    \draw[postaction={on each segment={mid arrow2}}]
      (p0) .. controls (p01c) .. (p1);
    \draw[red,postaction={on each segment={mid arrow3}}]
      (p1) .. controls (p12c) .. (p2);
    \draw[red,postaction={on each segment={mid arrow}}]
      (p3) .. controls (p32c) .. (p2);
  \end{scope}
  \filldraw[black] (p0) circle (.02*\tw)
    node[below] {\tiny $\pW_0$};
  \filldraw[black] (p1) circle (.02*\tw)
    node[below] {\tiny $\pW_1$};
  \filldraw[black] (p2) circle (.02*\tw)
    node[below] {\tiny $\pW_2$};
  \filldraw[black] (p3) circle (.02*\tw)
    node[below] {\tiny $\pW_3$};
\end{tikzpicture}}}
    \label{fig:torus-hyper4}
  }
\end{center}
\caption{(a) Identification of regions on~$\torus=[0,1]^2$ under the
  map~$\iota$.  (b) The surface~$\Sq=\torus/\iota$, with the four fixed points
  of~$\iota$ shown.  (c) $\Sq$ is a sphere with four punctures,
  denoted~$\S_{0,4}$.}
\end{figure}
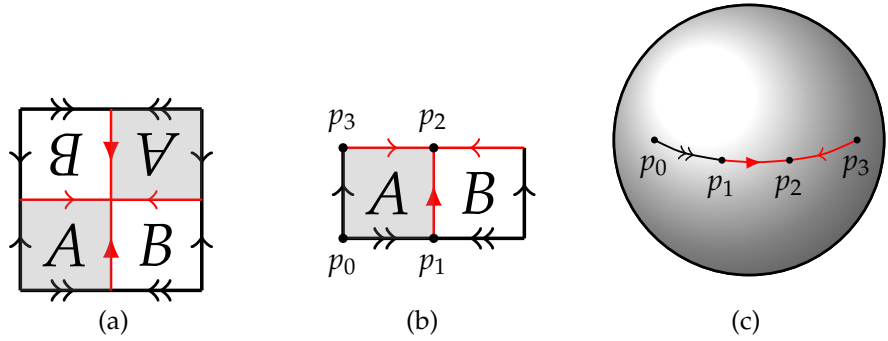

We now demonstrate that such a taffy puller motion arises naturally from
linear maps on the torus.  (We leave out many mathematical details --- see for
example \textcite{Fathi1979,Thurston1988,FarbMargalit} for the full story.)
We use the standard model of the torus~$\torus$ as the unit square~$[0,1]^2$
with opposite edges identified.  Consider the linear
map~$\iota:\torus\rightarrow\torus$, defined by~$\iota(x)=-x\mod 1$.  The
map~$\iota$ is an involution ($\iota^2=\id$) with four fixed points on the
torus~$[0,1]^2$,
\begin{equation}
  \pW_0 = \vect{0}{0}\,,\quad
  \pW_1 = \vect{\tfrac12}{0}\,,\quad
  \pW_2 = \vect{\tfrac12}{\tfrac12}\,,\quad
  \pW_3 = \vect{0}{\tfrac12}\,.
  \label{eq:pWdef}
\end{equation}
Figure~\ref{fig:torus-hyper} shows how the different sections of~$\torus$ are
mapped to each other under~$\iota$; arrows map to each other or are identified
because of periodicity.  The quotient space~$\Sq = \torus/\iota$, depicted in
Fig.~\ref{fig:torus-hyper2},
is actually a sphere in the topological sense (it has genus zero).  We can
see this by `gluing' the identified edges to obtain
Fig.~\ref{fig:torus-hyper4}.  The 4 fixed points of~$\iota$ above will play a
special role, so we puncture the sphere at those points and
write~$\Sq = \S_{0,4}$, which indicates a surface of genus~$0$ with~$4$
punctures.

Now let's take a general linear map~$\phi:\torus \rightarrow \torus$.  We
write~$\phi(x) = M\cdot x \mod 1$, with~$x \in [0,1]^2$ and~$M$ a matrix
in~$\SL_2(\Z)$,
\begin{equation}
  M = \begin{pmatrix}
    a & b \\ c & d
  \end{pmatrix},
  \qquad
  a,b,c,d \in \Z,\quad
  \quad ad-bc=1.
\end{equation}
This guarantees that~$\phi$ is an orientation-preserving homeomorphism --- a
continuous map of~$\torus$ with a continuous inverse.  The map~$\phi$
fixes~$\pW_0=\hvect{0}{0}$ and \emph{permutes} the ordered
set~$(\pW_1,\pW_2,\pW_3)$.  For example, the map
\begin{equation}
  \phi(x) = \begin{pmatrix}
    2 & 1 \\ 1 & 1
  \end{pmatrix}\cdot x\mod 1
  \label{eq:2111}
\end{equation}
maps~$(\pW_1,\pW_2,\pW_3)$ to~$(\pW_3,\pW_1,\pW_2)$.  This is an \emph{Anosov
  map}: it has a real eigenvalue larger than one (in magnitude).  We call the
spectral radius~$\lambda$ of the matrix~$M$ the \emph{dilatation} of the
map~$\phi$.  A key fact is that the length of any noncontractible simple
closed curve on~$\torus$ grows as~$\lambda^n$ as the number of
iterates~$n\rightarrow\infty$ \autocite{Fathi1979}.

Because~$\iota\circ\phi = \phi\circ\iota$, a linear map such as~$\phi$ on the
torus projects nicely to the punctured sphere~$\S_{0,4} = \torus/\iota$.  The
induced map on~$\S_{0,4}$ is called \emph{pseudo-Anosov} rather than Anosov,
since the quotient of the torus by~$\iota$ created four singularities.  (We
shall not need precise definitions of these terms; here by Anosov map we mean
a linear map on the torus with spectral radius larger than~$1$, and by
pseudo-Anosov we mean the same map projected to~$\S_{0,4}$.)

\begin{figure}
\begin{center}
  \subfigure[]{
    \resizebox{.65\subfigwidth}{!}{\begin{tikzpicture}
  \newcommand{\tw}{36pt}
  \begin{scope}[line width=.7pt]
    \draw (0,0) -- (\tw,0);
    \draw (0,\tw) -- (\tw,\tw);
    \draw (0,0) -- (0,\tw);
    \draw (\tw,0) -- (\tw,\tw);
  \end{scope}
  \begin{scope}[line width=.7pt]
    \draw[blue,thick,postaction={on each segment={mid arrow}}]
      (.5*\tw,.5*\tw) -- (0,.5*\tw);
    \draw[red,thick,postaction={on each segment={mid arrow}}]
      (.5*\tw,0) -- (.5*\tw,.5*\tw);

    \draw[black,fill=white] (0,0) circle (.03*\tw)
      node[below,scale=.75] {\tiny $\pW_0$};
    \filldraw[black] (1/2*\tw,0) circle (.02*\tw)
      node[below,scale=.75] {\tiny $\pW_1$};
    \filldraw[black] (0,1/2*\tw) circle (.02*\tw)
      node[left,scale=.75] {\tiny $\pW_3$};
    \filldraw[black] (1/2*\tw,1/2*\tw) circle (.02*\tw)
      node[right,scale=.75] {\tiny $\pW_2$};

    % Extra punctures using periodicity.
    \draw[black,fill=white] (1*\tw,0) circle (.03*\tw)
      node[below,scale=.75] {\tiny $\pW_0$};
    \filldraw[black] (1*\tw,1/2*\tw) circle (.02*\tw)
      node[right,scale=.75] {\tiny $\pW_3$};
    \draw[black,fill=white] (0,1*\tw) circle (.03*\tw)
      node[above,scale=.75] {\tiny $\pW_0$};
    \filldraw[black] (1/2*\tw,1*\tw) circle (.02*\tw)
      node[above,scale=.75] {\tiny $\pW_1$};
    \draw[black,fill=white] (1*\tw,1*\tw) circle (.03*\tw)
      node[above,scale=.75] {\tiny $\pW_0$};
  \end{scope}
\end{tikzpicture}}
  \hspace{2em}
    \raisebox{1em}{
    \resizebox{.8\subfigwidth}{!}{\begin{tikzpicture}
  \newcommand{\tw}{36pt}
  \coordinate (p0) at (0,0);
  \coordinate (p1) at (\tw,0);
  \coordinate (p2) at (2*\tw,0);
  \coordinate (p3) at (3*\tw,0);
  \begin{scope}[line width=1.2pt]
    % Draw outer ellipse, to suggest a sphere.
    \tikzbgellipsefour

    \draw[blue,ultra thick,postaction={on each segment={mid arrow}}]
      (p2) -- (p3);
    \draw[red,ultra thick,postaction={on each segment={mid arrow}}]
      (p1) -- (p2);

    \draw[black,fill=white] (p0) circle (.04*\tw)
      node[below] {\small $\pW_0$};
    \filldraw[black] (p1) circle (.03*\tw)
      node[below] {\small $\pW_1$};
    \filldraw[black] (p2) circle (.03*\tw)
      node[below] {\small $\pW_2$};
    \filldraw[black] (p3) circle (.03*\tw)
      node[below] {\small $\pW_3$};
  \end{scope}
\end{tikzpicture}}}
    \label{fig:torus-sphere-curve}
  }

  \subfigure[]{
    \resizebox{.65\subfigwidth}{!}{\begin{tikzpicture}
  \newcommand{\tw}{36pt}
  \begin{scope}[line width=.7pt]
    \draw (0,0) -- (\tw,0);
    \draw (0,\tw) -- (\tw,\tw);
    \draw (0,0) -- (0,\tw);
    \draw (\tw,0) -- (\tw,\tw);
  \end{scope}
  \begin{scope}[line width=.7pt]
    \draw[blue,thick,postaction={on each segment={mid arrow}}]
      (\tw,3/4*\tw) -- (1/2*\tw,1/2*\tw);
    \draw[blue,thick,postaction={on each segment={mid arrow}}]
      (1/2*\tw,\tw) -- (0,3/4*\tw);
    \draw[red,thick,postaction={on each segment={mid arrow}}]
      (0,1/2*\tw) -- (1/2*\tw,\tw);

    \draw[black,fill=white] (0,0) circle (.03*\tw)
      node[below,scale=.75] {\tiny $\pW_0'$};
    \filldraw[black] (1/2*\tw,0) circle (.02*\tw)
      node[below,scale=.75] {\tiny $\pW_2'$};
    \filldraw[black] (1/2*\tw,1/2*\tw) circle (.02*\tw)
      node[below right,scale=.75] {\tiny $\pW_3'$};
    \filldraw[black] (0,1/2*\tw) circle (.02*\tw)
      node[left,scale=.75] {\tiny $\pW_1'$};

    % Extra punctures using periodicity.
    \draw[black,fill=white] (1*\tw,0) circle (.03*\tw)
      node[below,scale=.75] {\tiny $\pW_0'$};
    \draw[black,fill=white] (0,1*\tw) circle (.03*\tw)
      node[above,scale=.75] {\tiny $\pW_0'$};
    \draw[black,fill=white] (1*\tw,1*\tw) circle (.03*\tw)
      node[above,scale=.75] {\tiny $\pW_0'$};
    \filldraw[black] (1/2*\tw,1*\tw) circle (.02*\tw)
      node[above,scale=.75] {\tiny $\pW_2'$};
    \filldraw[black] (1*\tw,1/2*\tw) circle (.02*\tw)
      node[right,scale=.75] {\tiny $\pW_1'$};
  \end{scope}
\end{tikzpicture}}
  \hspace{2em}
    \raisebox{1em}{
    \resizebox{.8\subfigwidth}{!}{%\usetikzlibrary{positioning}
\begin{tikzpicture}
  \newcommand{\tw}{36pt}
  \coordinate (p0) at (0,0);
  \coordinate (p1) at (\tw,0);
  \coordinate (p2) at (2*\tw,0);
  \coordinate (p3) at (3*\tw,0);
  \begin{scope}[line width=1.2pt]
    % Draw outer ellipse, to suggest a sphere.
    \tikzbgellipsefour

    %\draw[blue,ultra thick,postaction={on each segment={mid arrow}}]
    %  (p2) to[out=90,in=90] ($(p3) + (.2*\tw,0)$) to[out=-90,in=-90] (p1);
    \draw[blue,ultra thick,postaction={on each segment={mid arrow}}]
      (p1) to[out=-80,in=-90] ($(p3) + (.2*\tw,0)$) to[out=90,in=45] (p2);

    \draw[red,ultra thick,postaction={on each segment={mid arrow}}]
      (p3) to[out=-125,in=-35] (p1);

    \draw[black,fill=white] (p0) circle (.04*\tw)
      node[below = .05*\tw] {\small $\pW_0'$};
    \filldraw[black] (p1) circle (.03*\tw)
      node[below left] {\small $\pW_2'$};
    \filldraw[black] (p2) circle (.03*\tw)
      node[above left] {\small $\pW_3'$};
    \filldraw[black] (p3) circle (.03*\tw)
      node[left = .015*\tw] {\small $\pW_1'$};
  \end{scope}
\end{tikzpicture}}}
    \label{fig:torus-sphere-curve-2111}
  }
\end{center}
\caption{(a) Two curves on the torus~$\torus$ (left), which project to curves
  on the punctured sphere~$\S_{0,4}$ (right).  (b) The two curves transformed
  by the map~\eqref{eq:2111} (left), and projected onto~$\S_{0,4}$ (right).
  The transformed blue curve is the same as in the third frame of
  Fig.~\ref{fig:s1s-2_loop}.}
\label{fig:torus-sphere-2111}
\end{figure}
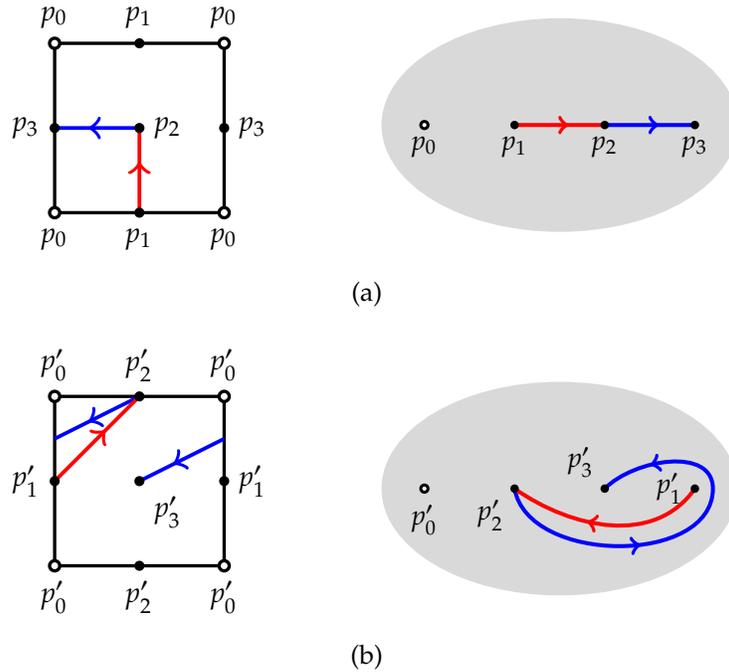

Let's see how the action of the map~\eqref{eq:2111} gives the taffy puller in
Fig.~\ref{fig:s1s-2}.  The permuted points~$(\pW_1,\pW_2,\pW_3)$ play the role
of the rods of the taffy puller.  Figure~\ref{fig:torus-sphere-curve} (left)
shows two curves on the torus, which project to curves on the punctured
sphere~$\S_{0,4}$ (right).  (Whenever we say curve, we will actually mean an
equivalence class of curves under homotopy fixing the punctures.)  The blue
curve from~$\pW_2$ to~$\pW_3$ should be identified with the piece of taffy in
Fig.~\ref{fig:s1s-2_loop}.  Now if we act on the curves with the torus
map~\eqref{eq:2111}, we obtain the curves in
Fig.~\ref{fig:torus-sphere-curve-2111} (left).  After taking the quotient
with~$\iota$, the curves project down as in
Fig.~\ref{fig:torus-sphere-curve-2111} (right).
This has the same shape as our taffy in Fig.~\ref{fig:s1s-2_loop} (third
frame) for $\tfrac13$ period of the taffy puller.

What we've essentially shown is that the taffy puller in Fig.~\ref{fig:s1s-2}
can be described by projecting the Anosov map~\eqref{eq:2111} of the torus to
a pseudo-Anosov map of~$\S_{0,4}$.  The growth of the length of taffy, under
repeated action, will be given by the spectral radius~$\lambda$ of the
matrix~$M$, here~$\lambda=\varphi^2$ with~$\varphi$ being the Golden Ratio
$\tfrac12(1+\sqrt5)$.  This taffy puller is a bit peculiar in that it requires
rods to move in a Figure-8 motion, as shown in
Fig.~\ref{fig:Nitz1918_rodmotion}.  This is challenging to achieve
mechanically, but surprisingly such a device was patented by
\textcite{nitz_candy-puller._1918} (Fig.~\ref{fig:Nitz1918_device}), and then
apparently again by \textcite{kirsch_candy-pulling_1928}.  The device requires
rods to alternately jump between two rotating wheels.

All 3-rod devices can be treated in the same manner, including the standard
3-rod taffy puller depicted in Fig.~\ref{fig:Robinson1908}.  We will not give
the details here, but it can be shown to arise from the linear map
\begin{equation}
  \phi(x) = \begin{pmatrix}
    5 & 2 \\ 2 & 1
  \end{pmatrix}\cdot x\mod 1
  \label{eq:5221}
\end{equation}
which has~$\lambda=\chi^2$.  Here~$\chi=1+\sqrt{2}$ is the \emph{Silver Ratio}
\autocite{MattFinn2011_silver}.

\section*{Are four rods better than three?}
%\label{sec:quad}

%
\begin{figure}
\begin{center}
\subfigure[]{
  \includegraphics[width=.5\textwidth]{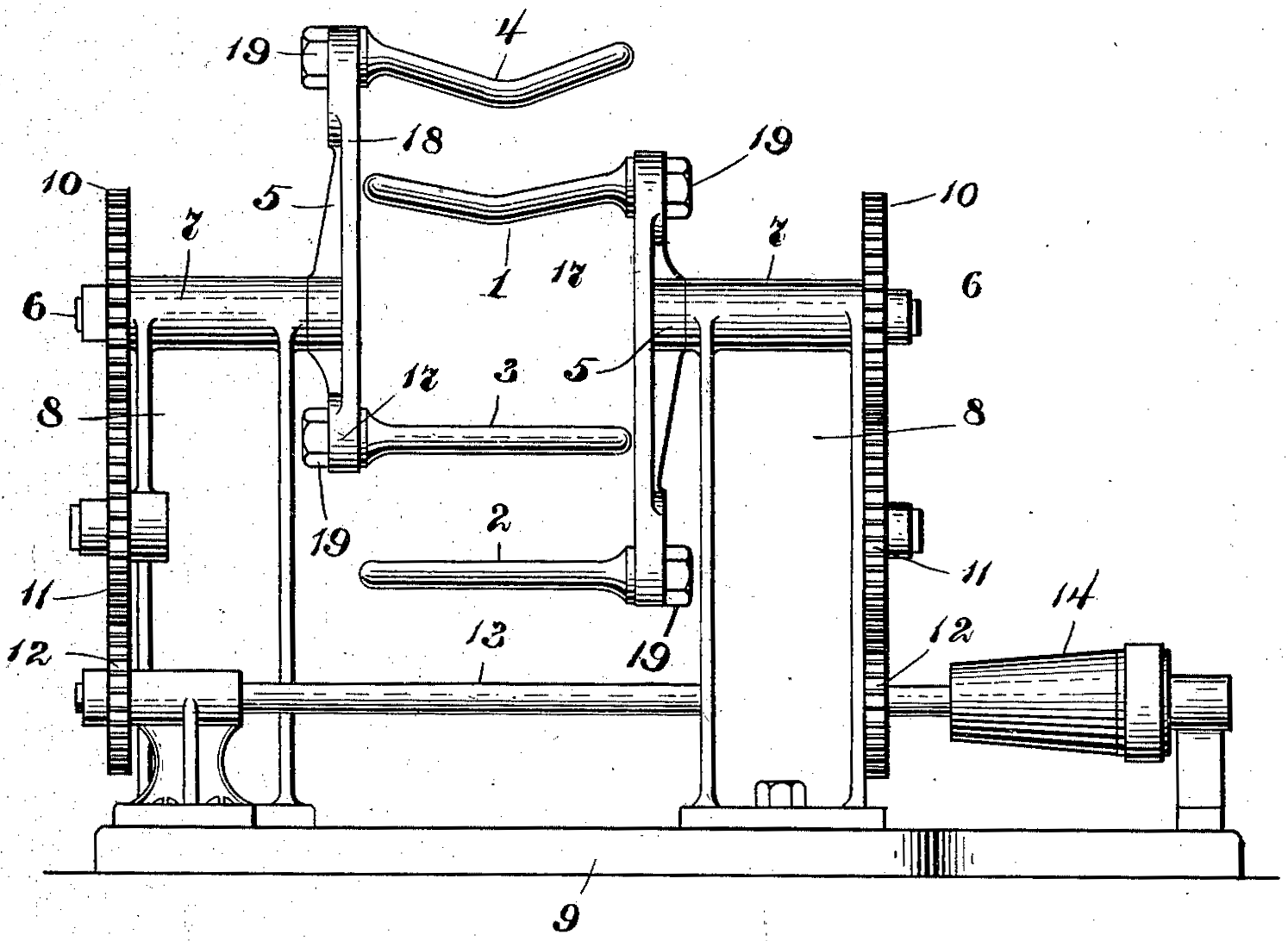}
  \label{fig:Thibodeau1903_device}
}%\hspace{1em}
\subfigure[]{
  \raisebox{1.7em}{
  \includegraphics[height=.2\textheight]{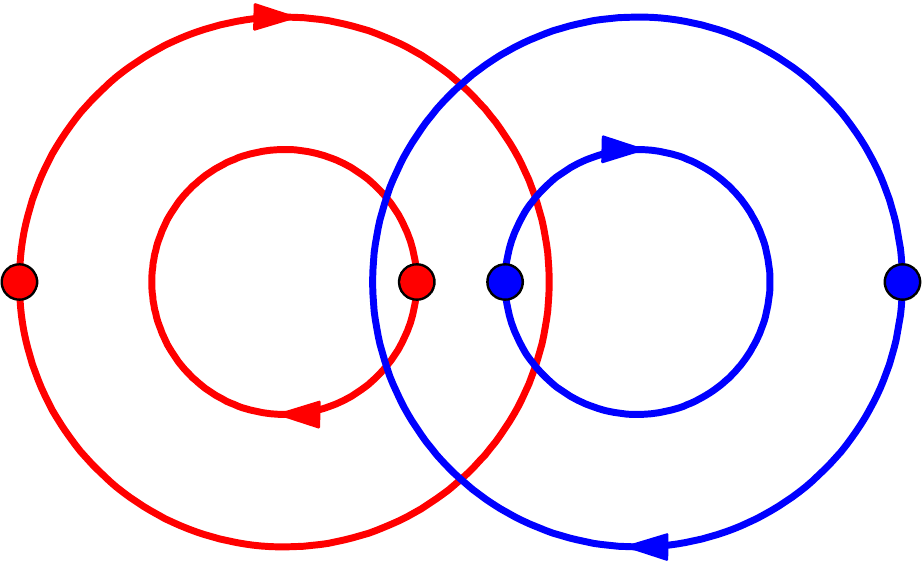}}
  \label{fig:Thibodeau1903_rodmotion}
}
\end{center}
\caption{(a) Side view of the standard 4-rod taffy puller from the patent
  of \textcite{thibodeau_method_1903}, with four rotating rods set on two
  axles.  (b) Rod motion.}
\label{fig:Thibodeau1903}
\end{figure}

As discussed in the previous section, all 3-rod taffy pullers arise from
Anosov maps of the torus.  This is not true in general for more than three
rods, but it is true for several specific devices.  Probably the most common
device is the \emph{standard 4-rod taffy puller}, which was invented by
\textcite{thibodeau_method_1903} and is shown in Fig.~\ref{fig:Thibodeau1903}.
It seems to have been rediscovered several times, such as by
\textcite{hudson_candy-working_1904}. \keepnote{There are more examples, but
  not sure it's worth mentioning.}  The design of
\textcite{richards_process_1905} is a variation that achieves the same effect,
and his patent has some of the prettiest diagrams of taffy pulling in action
(Fig.~\ref{fig:Richards1905}).  Mathematically, the 4-rod puller was studied
by \textcite{MacKay2001} and \textcite{Halbert2014}.

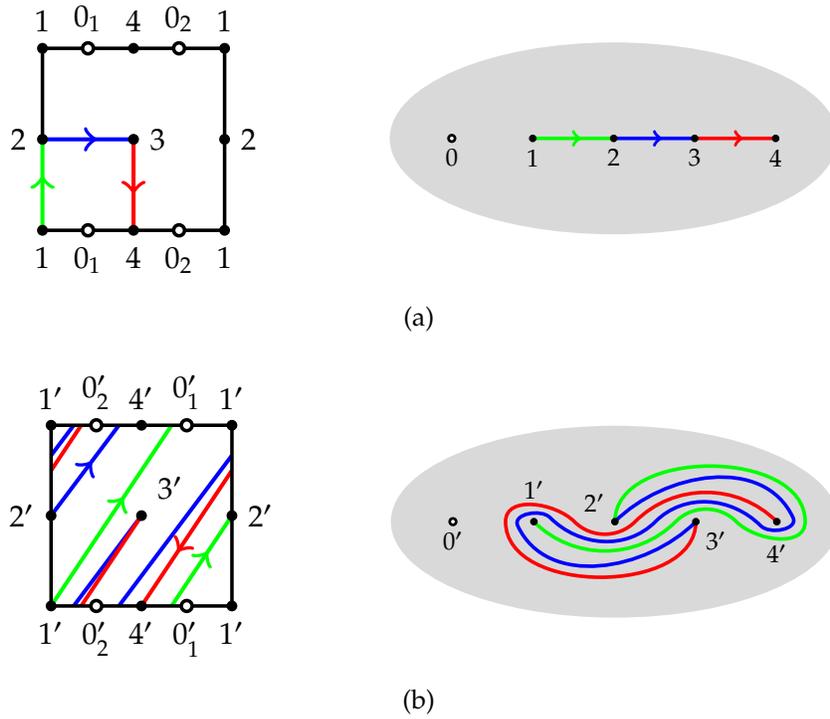
\begin{figure}
\begin{center}
  \subfigure[]{
    \resizebox{.64\subfigwidth}{!}{\begin{tikzpicture}
  \newcommand{\tw}{36pt}
  \coordinate (q01) at (1/4*\tw,0);
  \coordinate (q02) at (3/4*\tw,0);
  \begin{scope}[line width=.7pt]
    \draw[blue,thick,postaction={on each segment={mid arrow}}]
      (0,.5*\tw) -- (.5*\tw,.5*\tw);
    \draw[red,thick,postaction={on each segment={mid arrow}}]
      (.5*\tw,.5*\tw) -- (.5*\tw,0);
    \draw[green,thick,postaction={on each segment={mid arrow}}]
      (0,0) -- (0,.5*\tw);

    \begin{scope}[line width=.7pt]
      \draw (0,0) -- (\tw,0);
      \draw (0,\tw) -- (\tw,\tw);
      \draw (0,1/2*\tw) -- (0,\tw);
      \draw (\tw,0) -- (\tw,\tw);
    \end{scope}

    \filldraw[black] (0,0) circle (.02*\tw)
      node[below,scale=.75] {\tiny $1$};
    \filldraw[black] (1/2*\tw,0) circle (.02*\tw)
      node[below,scale=.75] {\tiny $4$};
    \filldraw[black] (0,1/2*\tw) circle (.02*\tw)
      node[left,scale=.75] {\tiny $2$};
    \filldraw[black] (1/2*\tw,1/2*\tw) circle (.02*\tw)
      node[right,scale=.75] {\tiny $3$};

    % Extra punctures using periodicity.
    \filldraw[black] (1*\tw,0) circle (.02*\tw)
      node[below,scale=.75] {\tiny $1$};
    \filldraw[black] (1*\tw,1/2*\tw) circle (.02*\tw)
      node[right,scale=.75] {\tiny $2$};
    \filldraw[black] (0,1*\tw) circle (.02*\tw)
      node[above,scale=.75] {\tiny $1$};
    \filldraw[black] (1/2*\tw,1*\tw) circle (.02*\tw)
      node[above,scale=.75] {\tiny $4$};
    \filldraw[black] (1*\tw,1*\tw) circle (.02*\tw)
      node[above,scale=.75] {\tiny $1$};

    % Boundary punctures.
    \draw[black,fill=white] (q01) circle (.03*\tw)
      node[below,scale=.75] {\tiny $0_1$};
    \draw[black,fill=white] (q02) circle (.03*\tw)
      node[below,scale=.75] {\tiny $0_2$};
    \draw[black,fill=white] ($(q01) + (0,\tw)$) circle (.03*\tw)
      node[above,scale=.75] {\tiny $0_1$};
    \draw[black,fill=white] ($(q02) + (0,\tw)$) circle (.03*\tw)
      node[above,scale=.75] {\tiny $0_2$};
  \end{scope}
\end{tikzpicture}}
  \hspace{2.5em}
    \raisebox{1.7em}{
    \resizebox{1\subfigwidth}{!}{\begin{tikzpicture}
  \newcommand{\tw}{36pt}
  \coordinate (q0) at (-1*\tw,0);
  \coordinate (q1) at (0,0);
  \coordinate (q2) at (\tw,0);
  \coordinate (q3) at (2*\tw,0);
  \coordinate (q4) at (3*\tw,0);
  \begin{scope}[line width=1.2pt]
    % Draw outer ellipse, to suggest a sphere.
    \tikzbgellipsefive

    \draw[green,ultra thick,postaction={on each segment={mid arrow}}]
      (q1) -- (q2);
    \draw[blue,ultra thick,postaction={on each segment={mid arrow}}]
      (q2) -- (q3);
    \draw[red,ultra thick,postaction={on each segment={mid arrow}}]
      (q3) -- (q4);

    \draw[black,fill=white] (q0) circle (.04*\tw)
      node[below] {\small $0$};
    \filldraw[black] (q1) circle (.03*\tw)
      node[below] {\small $1$};
    \filldraw[black] (q2) circle (.03*\tw)
      node[below] {\small $2$};
    \filldraw[black] (q3) circle (.03*\tw)
      node[below] {\small $3$};
    \filldraw[black] (q4) circle (.03*\tw)
      node[below] {\small $4$};
  \end{scope}
\end{tikzpicture}}
    }
    \label{fig:torus-sphere-curve2}
  }

  \subfigure[]{
    \resizebox{.68\subfigwidth}{!}{\begin{tikzpicture}
  \newcommand{\tw}{36pt}
  \coordinate (q01) at (1/4*\tw,0);
  \coordinate (q02) at (3/4*\tw,0);
  \begin{scope}[line width=.7pt]
    \draw[blue,thick,postaction={on each segment={mid arrow}}]
      (0,.5*\tw) -- (3/8*\tw,\tw);
    \draw[blue,thick]%,postaction={on each segment={mid arrow}}]
      (3/8*\tw,0) -- (\tw,5/6*\tw);
    \draw[blue,thick]%,postaction={on each segment={mid arrow}}]
      (0,5/6*\tw) -- (1/8*\tw,\tw);
    \draw[blue,thick]%,postaction={on each segment={mid arrow}}]
      (1/8*\tw,0) -- (1/2*\tw,1/2*\tw);

    \draw[red,thick]%,postaction={on each segment={mid arrow}}]
      (.5*\tw,.5*\tw) -- (1/6*\tw,0);
    \draw[red,thick]%,postaction={on each segment={mid arrow}}]
      (1/6*\tw,\tw) -- (0,3/4*\tw);
    \draw[red,thick,postaction={on each segment={mid arrow}}]
      (\tw,3/4*\tw) -- (1/2*\tw,0);

    \draw[green,thick,postaction={on each segment={mid arrow}}]
      (0,0) -- (2/3*\tw,\tw);
    \draw[green,thick,postaction={on each segment={mid arrow}}]
      (2/3*\tw,0) -- (\tw,1/2*\tw);

    \begin{scope}[line width=.7pt]
      \draw (0,0) -- (\tw,0);
      \draw (0,\tw) -- (\tw,\tw);
      \draw (0,0) -- (0,\tw);
      \draw (\tw,0) -- (\tw,\tw);
    \end{scope}

    \filldraw[black] (0,0) circle (.02*\tw)
      node[below,scale=.75] {\tiny $1'$};
    \filldraw[black] (1/2*\tw,0) circle (.02*\tw)
      node[below,scale=.75] {\tiny $4'$};
    \filldraw[black] (0,1/2*\tw) circle (.02*\tw)
      node[left,scale=.75] {\tiny $2'$};
    \filldraw[black] (1/2*\tw,1/2*\tw) circle (.02*\tw)
      node[above right,scale=.75] {\tiny $3'$};

    % Extra punctures using periodicity.
    \filldraw[black] (1*\tw,0) circle (.02*\tw)
      node[below,scale=.75] {\tiny $1'$};
    \filldraw[black] (1*\tw,1/2*\tw) circle (.02*\tw)
      node[right,scale=.75] {\tiny $2'$};
    \filldraw[black] (0,1*\tw) circle (.02*\tw)
      node[above,scale=.75] {\tiny $1'$};
    \filldraw[black] (1/2*\tw,1*\tw) circle (.02*\tw)
      node[above,scale=.75] {\tiny $4'$};
    \filldraw[black] (1*\tw,1*\tw) circle (.02*\tw)
      node[above,scale=.75] {\tiny $1'$};

    % Boundary punctures.
    \draw[black,fill=white] (q01) circle (.03*\tw)
      node[below,scale=.75] {\tiny $0_2'$};
    \draw[black,fill=white] (q02) circle (.03*\tw)
      node[below,scale=.75] {\tiny $0_1'$};
    \draw[black,fill=white] ($(q01) + (0,\tw)$) circle (.03*\tw)
      node[above,scale=.75] {\tiny $0_2'$};
    \draw[black,fill=white] ($(q02) + (0,\tw)$) circle (.03*\tw)
      node[above,scale=.75] {\tiny $0_1'$};
  \end{scope}
\end{tikzpicture}}
  \hspace{2em}
    \raisebox{1.7em}{
    \resizebox{1\subfigwidth}{!}{\begin{tikzpicture}
  \newcommand{\tw}{36pt}
  \coordinate (q0) at (-1*\tw,0);
  \coordinate (q1) at (0,0);
  \coordinate (q2) at (\tw,0);
  \coordinate (q3) at (2*\tw,0);
  \coordinate (q4) at (3*\tw,0);
  \begin{scope}[line width=1.2pt]
    % Draw outer ellipse, to suggest a sphere.
    \tikzbgellipsefive

    \draw[green,ultra thick]%,postaction={on each segment={mid arrow}}]
      (q1)
      to[out=-45,in=-135] ($(q2) + (.75*\tw,0)$)
      to[out=45,in=135] ($(q3) + (.5*\tw,0)$)
      to[out=-45,in=-90] ($(q4) + (.35*\tw,0)$)
      to[out=90,in=90] (q2)
      ;

    \draw[red,ultra thick]%,postaction={on each segment={mid arrow}}]
      (q4)
      to[out=135,in=45] ($(q2) + (.25*\tw,0)$)
      to[out=-135,in=-45] ($(q1) + (.5*\tw,0)$)
      to[out=135,in=90] ($(q1) - (.35*\tw,0)$)
      to[out=-90,in=-90] (q3)
      ;

    \draw[blue,ultra thick]%,postaction={on each segment={mid arrow}}]
      (q2)
      to[out=45,in=115] ($(q4) + (.2*\tw,.05*\tw)$)
      to[out=-65,in=-45] ($(q4) - (.2*\tw,.05*\tw)$)
      to[out=135,in=45] ($(q2) + (.5*\tw,0)$)
      to[out=-135,in=-45] ($(q1) + (.2*\tw,.05*\tw)$)
      to[out=135,in=115] ($(q1) - (.2*\tw,.05*\tw)$)
      to[out=-65,in=-135] (q3)
      ;

    \draw[black,fill=white] (q0) circle (.04*\tw)
      node[below] {\small $0'$};
    \filldraw[black] (q1) circle (.03*\tw)
      node[above = .15*\tw] {\small $1'$};
    \filldraw[black] (q2) circle (.03*\tw)
      node[above left] {\small $2'$};
    \filldraw[black] (q3) circle (.03*\tw)
      node[below right] {\small $3'$};
    \filldraw[black] (q4) circle (.03*\tw)
      node[below = .15*\tw] {\small $4'$};
  \end{scope}
\end{tikzpicture}}
    }
    \label{fig:torus-sphere-curve2-3243}
  }
\end{center}
\caption{(a) Three curves on the torus~$\torus$ (left), which project to
  curves on the punctured sphere~$\S_{0,5}$ (right).  (b) The three curves
  transformed by the map~\eqref{eq:3243} (left), and projected onto~$\S_{0,5}$
  (right).  Compare these to the last frame of Fig.~\ref{fig:Richards1905}.}
\end{figure}

The rod motion for the standard 4-rod puller is shown in
Fig.~\ref{fig:Thibodeau1903_rodmotion}.  Observe that the two orbits of
smaller radius are not intertwined, so topologically they might as well be
fixed rods.  \keepnote{A full period of the puller requires that all the rods
  to return to their initial position.} %
This taffy puller arises from an Anosov map such as~\eqref{eq:5221}, but with
all four points~$(\pW_0,\pW_1,\pW_2,\pW_3)$ of~$\S_{0,4}$ identified with
rods.  We relabel the four points~$(\pW_0,\pW_1,\pW_2,\pW_3)$ as~$(1,4,3,2)$,
as in Fig.~\ref{fig:torus-sphere-curve2} (left), which gives the order of the
rods on the right in that figure.  The boundary point labeled~$0$ plays no
direct role, but prevents us from shortening curves by passing them `behind'
the sphere.%
\footnote{%
  In the 3-rod case, the point labeled~$0$ in Fig.~\ref{fig:torus-sphere-2111}
  plays this role.  In the 4-rod case, we need to use~$\S_{0,5}$ in order to
  allow for this extra point.  There are no more fixed points available,
  since~$\phi(x)$ in~\eqref{eq:3243} only has 4.  However, a period-2 point
  of~$\phi$ will do, as long as the two iterates are also mapped to each other
  by~$\iota$.  The map~\eqref{eq:3243} actually has 14 orbits of period 2, but
  only two of those are also invariant under~$\iota$:
  \begin{equation*}
    \l\{\vect{\tfrac14}{0}\,,
        \vect{\tfrac34}{0} \r\}
    \quad\text{and}\quad
    \l\{\vect{\tfrac14}{\tfrac12}\,,
        \vect{\tfrac34}{\tfrac12} \r\}.
  \end{equation*}
  The second choice would put the boundary point being between two rods, so we
  choose the first orbit.  The two iterates are labeled~$0_1$ and~$0_2$ in
  Fig.~\ref{fig:torus-sphere-curve2}.  They are interchanged in
  Fig.~\ref{fig:torus-sphere-curve2-3243} after applying the
  map~\eqref{eq:3243}, but they both map to the same point on the
  sphere~$\S_{0,5} = (\torus/\iota) - \{0\}$, since they also
  satisfy~$\iota(0_1)=0_2$. %
} %
Puncturing at this extra point gives the space~$\S_{0,5}$, the sphere with~$5$
punctures.

Now act on the curves in Fig.~\ref{fig:torus-sphere-curve2} with the map
\begin{equation}
  \phi(x) = \begin{pmatrix}
    3 & 2 \\ 4 & 3
  \end{pmatrix}\cdot x\mod 1.
  \label{eq:3243}
\end{equation}
This map fixes each of the points~$1,2,3,4$, just as the 4-rod taffy puller
does.  Figure~\ref{fig:torus-sphere-curve2-3243} shows the action of the map
on curves anchored on the rods: it acts in exactly the same manner as the
standard 4-rod taffy puller.  In fact,
\begin{equation}
  \begin{pmatrix}
    3 & 2 \\ 4 & 3
  \end{pmatrix}
  =
  \begin{pmatrix} 1 & 0 \\ 1 & 1 \end{pmatrix}
  \begin{pmatrix}
    5 & 2 \\ 2 & 1
  \end{pmatrix}
  \begin{pmatrix} 1 & 0 \\ 1 & 1 \end{pmatrix}^{-1}
\end{equation}
which means that the maps~\eqref{eq:5221} and~\eqref{eq:3243} are
\emph{conjugate} to each other.  Conjugate maps have the same dilatation (the
trace is invariant), so the standard 3-rod and 4-rod taffy pullers arise from
essentially the same Anosov map, only interpreted differently.  In other
words, at least for the standard 4-rod puller, the addition of a rod does
\emph{not} increase the effectiveness of the device.

\section*{A new device}

All the devices we described so far arise from maps of the torus.  Now we give
an example of a device that arises from a \emph{branched cover} of the torus,
rather than directly from the torus itself.  (A theorem of
\textcite{Franks1999} implies that the dilatation~$\lambda$ must also be
quadratic in this case.)
\begin{figure}
\begin{center}
\subfigure[]{
  \includegraphics[width=.45\textwidth]{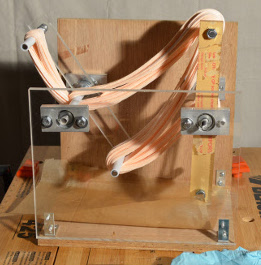}
  \label{fig:6rods_photo}
}
\hspace{1em}
\subfigure[]{
  \raisebox{3em}{%
    \includegraphics[width=.45\textwidth]{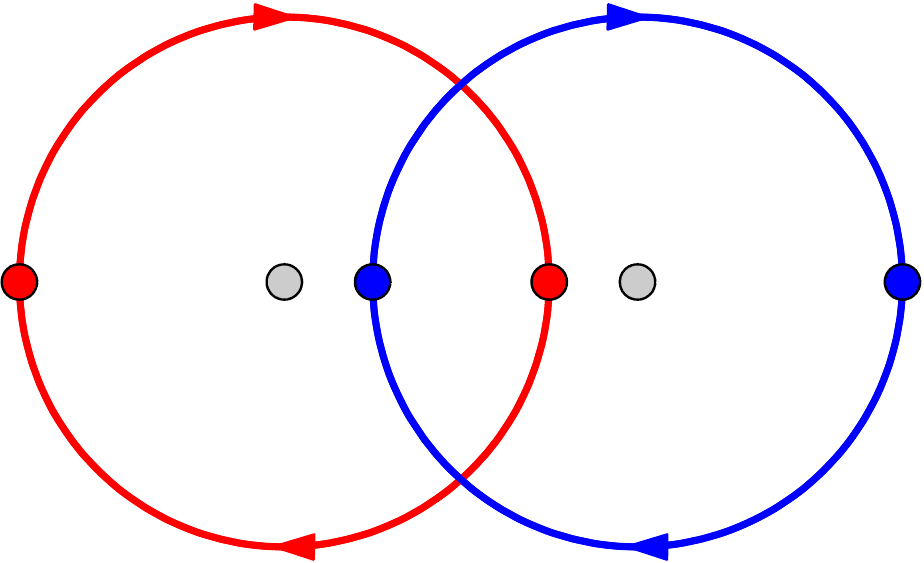}}
  \label{fig:6rods_rodmotion}
}
\end{center}
\caption{(a) A 6-pronged taffy puller designed and built by Alexander
  Flanagan and the author.  (b) The motion of the rods, with two fixed axles.}
\label{fig:6rods}
\end{figure}
Figure~\ref{fig:6rods} shows such a device, designed and built by Alexander
Flanagan and the author.  It is a simple modification of the standard 4-rod
design (Fig.~\ref{fig:Thibodeau1903}), except that the two arms are of equal
length, and the axles are extended to become fixed rods.  There are thus 6
rods in play, and we shall see that this device has a rather large dilatation.

\begin{figure}
\newlength{\invfigwidth}
\setlength{\invfigwidth}{1.025\textwidth}
\begin{center}
  \subfigure[]{
    \includegraphics[width=.4\invfigwidth]{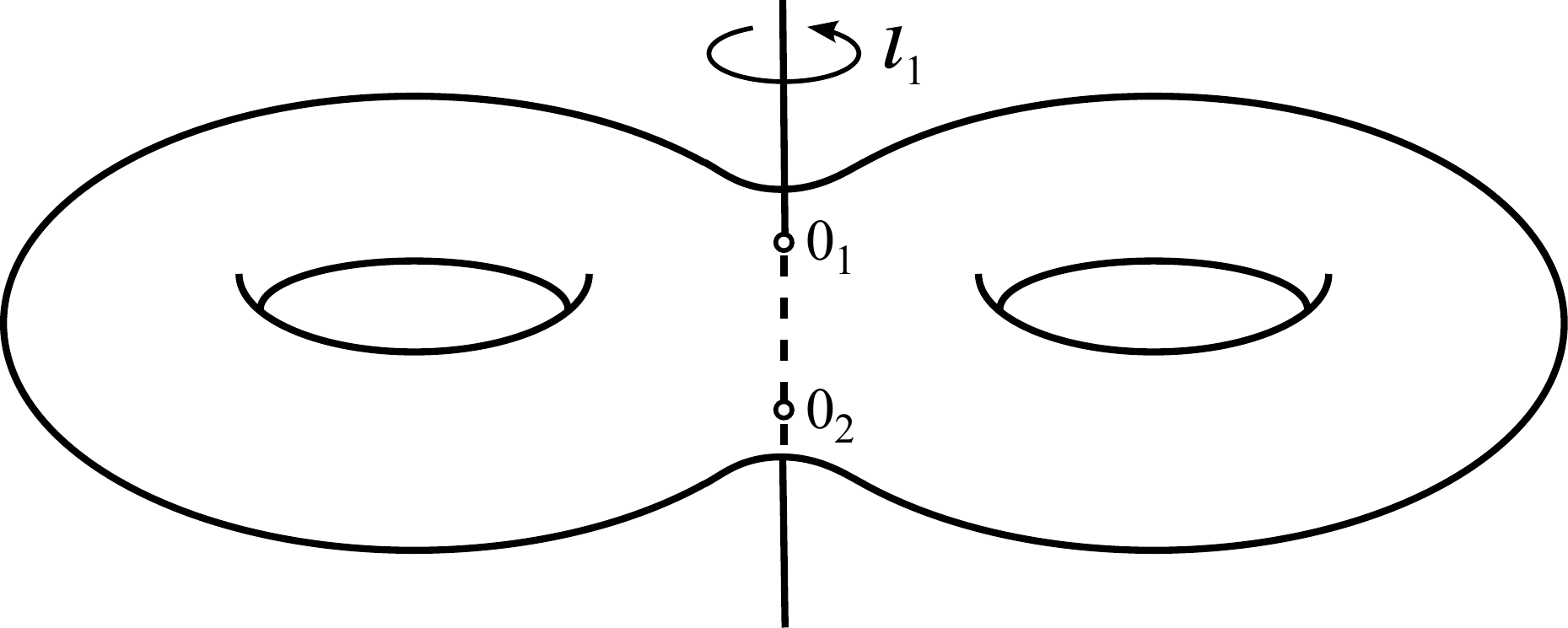}
    \label{fig:2torus_inv1}
  }
  \hspace{.25em}
  \subfigure[]{
    \raisebox{.018\invfigwidth}{
      \includegraphics[width=.51\invfigwidth]{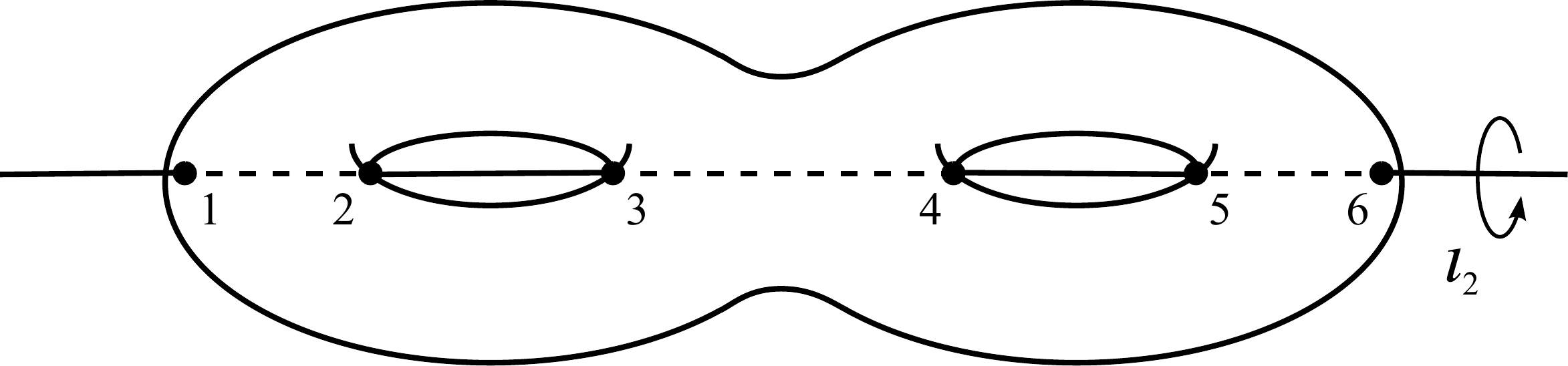}
    }
    \label{fig:2torus_inv2}
  }
\end{center}
\caption{The two involutions of a genus two surface~$\S_2$ as rotations
  by~$\pi$.  (a) The involution~$\iota_1$ has two fixed points; (b) $\iota_2$
  has six.}
\label{fig:genus2_involutions}
\end{figure}

The construction of a map describing this 6-rod device uses the two
involutions of the closed (unpunctured) genus two surface~$\S_2$ shown in
Fig.~\ref{fig:genus2_involutions}.  Imagine that an Anosov map gives the
dynamics on the left `torus' of the surface.  The involution~$\iota_1$ extends
those dynamics to a genus two surface.  The involution~$\iota_2$ is then used
to create the quotient surface~$\S_{0,6}=\S_2/\iota_2$.  The 6 punctures will
correspond to the rods of the taffy puller.

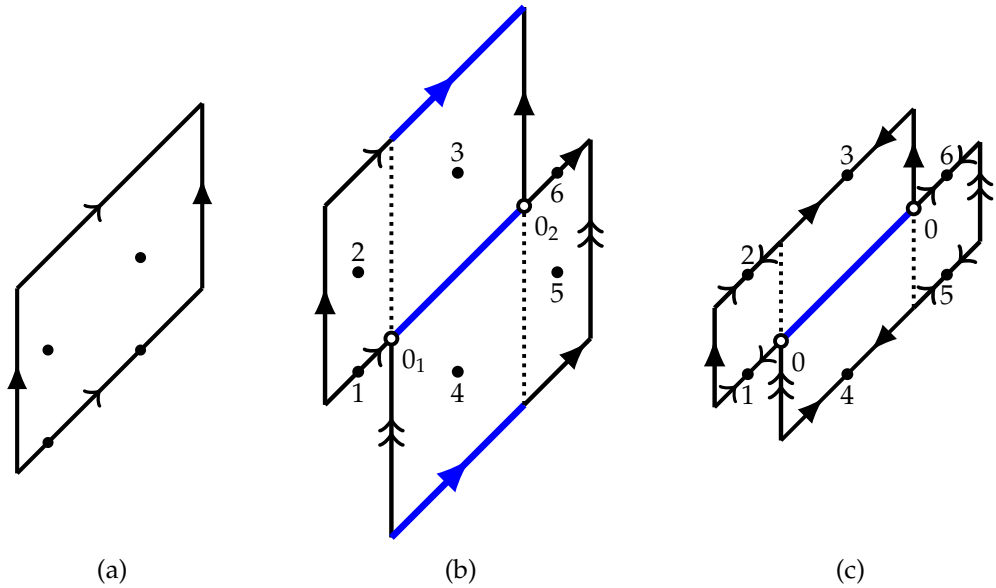
\begin{figure}
\newlength{\ssubfigwidth}
\setlength{\ssubfigwidth}{.9\subfigwidth}
\begin{center}
  \subfigure[]{
    \raisebox{.16\ssubfigwidth}{
      \resizebox{.5\ssubfigwidth}{!}{\begin{tikzpicture}
  \newcommand{\tw}{36pt}

  \coordinate (O) at (0,0);
  \coordinate (p1) at (-1/2*\tw,-1/2*\tw);
  \coordinate (p2) at (-1/2*\tw,0);
  \coordinate (p3) at (0,1/2*\tw);

  \coordinate (ll) at (-2/3*\tw,-2/3*\tw);
  \coordinate (ul) at (-2/3*\tw, 1/3*\tw);
  \coordinate (lr) at ( 1/3*\tw, 1/3*\tw);
  \coordinate (ur) at ( 1/3*\tw, 4/3*\tw);

  \begin{scope}[line width=.7pt]

    % Top torus
    \draw[thick,-|>-=.6] (ll) -- (ul);
    \draw[thick,-|>-=.6] (lr) -- (ur);
    \draw[thick,->-=.45] (ul) -- (ur);
    \draw[thick,->-=.45] (ll) -- (lr);

    \filldraw[black] (O) circle (.02*\tw);
      %node[below,scale=.65] {\tiny $\,\,\,\,\,\,\,(0,0)$};
    \filldraw[black] (p1) circle (.02*\tw);
      %node[below,scale=.75] {\tiny $2$};
    \filldraw[black] (p2) circle (.02*\tw);
      %node[above,scale=.75] {\tiny $1$};
    \filldraw[black] (p3) circle (.02*\tw);
      %node[above,scale=.75] {\tiny $3$};

  \end{scope}
\end{tikzpicture}}
    }
    \label{fig:genus2_half}
  }
  \hspace{2em}
  \subfigure[]{
    \resizebox{.7\ssubfigwidth}{!}{\begin{tikzpicture}
  \newcommand{\tw}{36pt}

  \coordinate (O) at (0,0);
  \coordinate (p1) at (-1/2*\tw,-1/2*\tw);
  \coordinate (p2) at (-1/2*\tw,0);
  \coordinate (p3) at (0,1/2*\tw);
  \coordinate (p4) at (0,-1/2*\tw);
  \coordinate (p5) at (1/2*\tw,0);
  \coordinate (p6) at (1/2*\tw,1/2*\tw);

  \coordinate (ll) at (-2/3*\tw,-2/3*\tw);
  \coordinate (ul) at (-2/3*\tw, 1/3*\tw);
  \coordinate (lr) at ( 1/3*\tw, 1/3*\tw);
  \coordinate (ur) at ( 1/3*\tw, 4/3*\tw);

  \begin{scope}[line width=.7pt]

    % Top torus
    \draw[thick,-|>-=.6] (ll) -- (ul);
    \draw[thick,-|>-=.6] (lr) -- (ur);
    \draw[thick,->-=.85] (ul) -- ($(ul) + (1/3*\tw,1/3*\tw)$);
    \draw[thick,->-=.85] (ll) -- ($(ll) + (1/3*\tw,1/3*\tw)$);
    \draw[blue, very thick,-|>-=.5] ($(ul) + (1/3*\tw,1/3*\tw)$) -- (ur);
    \draw[blue, very thick]  ($(ll) + (1/3*\tw,1/3*\tw)$) -- (lr);
    \draw[densely dotted]
      ($(ll) + (1/3*\tw,1/3*\tw)$) -- ($(ul) + (1/3*\tw,1/3*\tw)$);

    % Bottom torus
    \coordinate (sh) at (1/3*\tw,-2/3*\tw);
    \draw[thick,->>-=.6] ($(ll)+(sh)$) -- ($(ul)+(sh)$);
    \draw[thick,->>-=.6] ($(lr)+(sh)$) -- ($(ur)+(sh)$);
    \draw[thick,-|>-=.95] ($(ur)+(sh) - (1/3*\tw,1/3*\tw)$) -- ($(ur)+(sh)$);
    \draw[thick,-|>-=.95] ($(lr)+(sh) - (1/3*\tw,1/3*\tw)$) -- ($(lr)+(sh)$);
    \draw[blue, very thick,-|>-=.6]
      ($(ll)+(sh)$) -- ($(ll)+(sh) + (2/3*\tw,2/3*\tw)$);
    \draw[densely dotted]
      ($(lr)+(sh) - (1/3*\tw,1/3*\tw)$) -- ($(ur)+(sh) - (1/3*\tw,1/3*\tw)$);

    %\draw[black,fill=white] (O) circle (.03*\tw)
    %  node[below,scale=.65] {\tiny $\,\,\,\,\,\,\,(0,0)$};
    \draw[black,fill=white] ($(O) - (1/3*\tw,1/3*\tw)$) circle (.03*\tw)
      node[below right,scale=.45] {\small $0_1$};
    \draw[black,fill=white] ($(O) + (1/3*\tw,1/3*\tw)$) circle (.03*\tw)
      node[below right,scale=.45] {\small $0_2$};
    \filldraw[black] (p1) circle (.02*\tw)
      node[below,scale=.45] {\small $1$};
    \filldraw[black] (p2) circle (.02*\tw)
      node[above,scale=.45] {\small $2$};
    \filldraw[black] (p3) circle (.02*\tw)
      node[above,scale=.45] {\small $3$};
    \filldraw[black] (p4) circle (.02*\tw)
      node[below,scale=.45] {\small $4$};
    \filldraw[black] (p5) circle (.02*\tw)
      node[below,scale=.45] {\small $5$};
    \filldraw[black] (p6) circle (.02*\tw)
      node[below,scale=.45] {\small $6$};

  \end{scope}
\end{tikzpicture}}
    \label{fig:genus2_flat}
  }
  \hspace{2em}
  \subfigure[]{
    \raisebox{.24\ssubfigwidth}{
      \resizebox{.7\ssubfigwidth}{!}{\begin{tikzpicture}
  \newcommand{\tw}{36pt}

  \coordinate (O) at (0,0);
  \coordinate (p1) at (-1/2*\tw,-1/2*\tw);
  \coordinate (p2) at (-1/2*\tw,0);
  \coordinate (p3) at (0,1/2*\tw);
  \coordinate (p4) at (0,-1/2*\tw);
  \coordinate (p5) at (1/2*\tw,0);
  \coordinate (p6) at (1/2*\tw,1/2*\tw);

  \coordinate (ll) at (-2/3*\tw,-2/3*\tw);
  \coordinate (ul) at (-2/3*\tw,-1/6*\tw);
  \coordinate (lr) at ( 1/3*\tw, 1/3*\tw);
  \coordinate (ur) at ( 1/3*\tw, 5/6*\tw);

  \begin{scope}[line width=.7pt]

    % Top torus
    % Periodic edges
    \draw[thick,-|>-=.7] (ll) -- (ul);
    \draw[thick,-|>-=.7] (lr) -- (ur);
    % Puncture 1
    \draw[thick,->-=.65] (ul) -- ($(ul) + 1/2*(1/3*\tw,1/3*\tw)$);
    \draw[thick,->-=.65]
      ($(ul) + (1/3*\tw,1/3*\tw)$) -- ($(ul) + 1/2*(1/3*\tw,1/3*\tw)$);
    % Puncture 2
    \draw[->-=.65] (ll) -- ($(ll) + 1/2*(1/3*\tw,1/3*\tw)$);
    \draw[->-=.65]
      ($(ll) + (1/3*\tw,1/3*\tw)$) -- ($(ll) + 1/2*(1/3*\tw,1/3*\tw)$);
    % Puncture 3
    \draw[-|>-=.65]
      ($(ul) + (1/3*\tw,1/3*\tw)$) -- ($(ul) + (2/3*\tw,2/3*\tw)$);
    \draw[-|>-=.65]
      (ur) -- ($(ul) + (2/3*\tw,2/3*\tw)$);

    \draw[blue, very thick]  ($(ll) + (1/3*\tw,1/3*\tw)$) -- (lr);
    \draw[densely dotted]
      ($(ll) + (1/3*\tw,1/3*\tw)$) -- ($(ul) + (1/3*\tw,1/3*\tw)$);

    % Bottom torus
    \coordinate (sh) at (1/3*\tw,-1/6*\tw);
    % Periodic edges
    \draw[thick,->>-=.7] ($(ll)+(sh)$) -- ($(ul)+(sh)$);
    \draw[thick,->>-=.7] ($(lr)+(sh)$) -- ($(ur)+(sh)$);
    % Puncture 5
    \draw[thick,->-=.65] ($(ur)+(sh)$) -- ($(ur)+(sh) - 1/2*(1/3*\tw,1/3*\tw)$);
    \draw[thick,->-=.65]
      ($(ur)+(sh) - (1/3*\tw,1/3*\tw)$) -- ($(ur)+(sh) - 1/2*(1/3*\tw,1/3*\tw)$);   % Puncture 6
    \draw[thick,->-=.65] ($(lr)+(sh)$) -- ($(lr)+(sh) - 1/2*(1/3*\tw,1/3*\tw)$);
    \draw[thick,->-=.65]
      ($(lr)+(sh) - (1/3*\tw,1/3*\tw)$) -- ($(lr)+(sh) - 1/2*(1/3*\tw,1/3*\tw)$);
    % Puncture 4
    \draw[-|>-=.65]
      ($(ll)+(sh)$) -- ($(ll)+(sh) + (1/3*\tw,1/3*\tw)$);
    \draw[-|>-=.65]
      ($(ll)+(sh) + (2/3*\tw,2/3*\tw)$) -- ($(ll)+(sh) + (1/3*\tw,1/3*\tw)$);

    % Bottom torus
    %\draw[thick,-|>-=.95] ($(ur)+(sh) - (1/3*\tw,1/3*\tw)$) -- ($(ur)+(sh)$);
    %\draw[thick,-|>-=.95] ($(lr)+(sh) - (1/3*\tw,1/3*\tw)$) -- ($(lr)+(sh)$);
    %\draw[blue, very thick,-|>-=.6]
     % ($(ll)+(sh)$) -- ($(ll)+(sh) + (2/3*\tw,2/3*\tw)$);
    \draw[densely dotted]
      ($(lr)+(sh) - (1/3*\tw,1/3*\tw)$) -- ($(ur)+(sh) - (1/3*\tw,1/3*\tw)$);

    %\draw[black,fill=white] (O) circle (.03*\tw)
    %  node[below,scale=.65] {\tiny $\,\,\,\,\,\,\,(0,0)$};
    \draw[black,fill=white] ($(O) - (1/3*\tw,1/3*\tw)$) circle (.03*\tw)
      node[below right,scale=.45] {\small $0$}; % was 0_1
    \draw[black,fill=white] ($(O) + (1/3*\tw,1/3*\tw)$) circle (.03*\tw)
      node[below right,scale=.45] {\small $0$}; % was 0_2
    \filldraw[black] (p1) circle (.02*\tw)
      node[below,scale=.45] {\small $1$};
    \filldraw[black] (p2) circle (.02*\tw)
      node[above,scale=.45] {\small $2$};
    \filldraw[black] (p3) circle (.02*\tw)
      node[above,scale=.45] {\small $3$};
    \filldraw[black] (p4) circle (.02*\tw)
      node[below,scale=.45] {\small $4$};
    \filldraw[black] (p5) circle (.02*\tw)
      node[below,scale=.45] {\small $5$};
    \filldraw[black] (p6) circle (.02*\tw)
      node[above,scale=.45] {\small $6$};

  \end{scope}
\end{tikzpicture}}
    }
    \label{fig:genus2_sphere}
  }
\end{center}
\caption{(a) A different unfolding of the torus.  The four fixed points
  of~$\iota$ are indicated.  (b) Two copies of the torus glued together after
  removing a disk.  The points~$0_{1,2}$ are
  at~$\mp\hvect{\tfrac13}{\tfrac13}$.  This gives the genus two
  surface~$\S_2$.  The two tori are mapped to each other by the
  involution~$\iota_1$ from Fig~\ref{fig:genus2_involutions}, with fixed
  points~$0_{1,2}$. The involution~$\iota_2$ acts on the individual tori with
  fixed points~$1,\ldots,6$.  (c) The quotient surface~$\S_2/\iota_2$, which
  is the punctured sphere~$\S_{0,6}$.}
%\label{fig:genus2_sphere_flat}
\end{figure}

A bit of experimentation suggests starting from the Anosov map
\begin{equation}
  \phi(x) = \begin{pmatrix}
    -1 & -1 \\ -2 & -3
  \end{pmatrix}\cdot x\mod 1.
  \label{eq:1123}
\end{equation}
Referring to the points~\eqref{eq:pWdef}, this map fixes~$\pW_0$ and~$\pW_1$
and interchanges~$\pW_2$ and~$\pW_3$.  For our purposes, we cut our unit cell
for the torus slightly differently, as shown in Fig.~\ref{fig:genus2_half}.
In addition to~$\pW_0$ and~$\pW_1$, the map has four more fixed points:
\begin{equation}
  %\l\{
  \vect{-\tfrac13}{-\tfrac13}\,,\,
  \vect{\tfrac13}{\tfrac13}\,,\,
  \vect{\tfrac16}{\tfrac23}\,,\,
  \vect{-\tfrac16}{\phantom{-}\tfrac13}
  %\r\}
  \,.
\end{equation}
To create our branched cover of the torus, we will make a cut from the
point~$0_1 = \hvect{-\tfrac13}{-\tfrac13}$
to~$0_2 = \hvect{\tfrac13}{\tfrac13}$, as shown in Fig.~\ref{fig:genus2_flat}.
We have also labeled by~$1$--$6$ the points that will correspond to our rods.
The arrows show identified opposite edges; we have effectively cut a slit in
two tori, opened the slits into disks, and glued the tori at those disks to
create a genus two surface.  The involution~$\iota_1$ from
Fig.~\ref{fig:2torus_inv1} corresponds to translating the top half in
Fig.~\ref{fig:genus2_flat} down to the bottom half; the only fixed points are
then~$0_1$ and~$0_2$.  For the involution~$\iota_2$ of
Fig.~\ref{fig:2torus_inv2}, first divide Fig.~\ref{fig:genus2_flat} into four
sectors with~$2$--$5$ at their center; then rotate each sector by~$\pi$ about
its center.  This fixes the points~$1$--$6$.  The quotient
surface~$\S_2/\iota_2$ gives the punctured sphere~$\S_{0,6}$, shown in
Fig.~\ref{fig:genus2_sphere}.  The points~$0_{1,2}$ are mapped to each other
by~$\iota_2$ and so become identified with the same point~$0$.

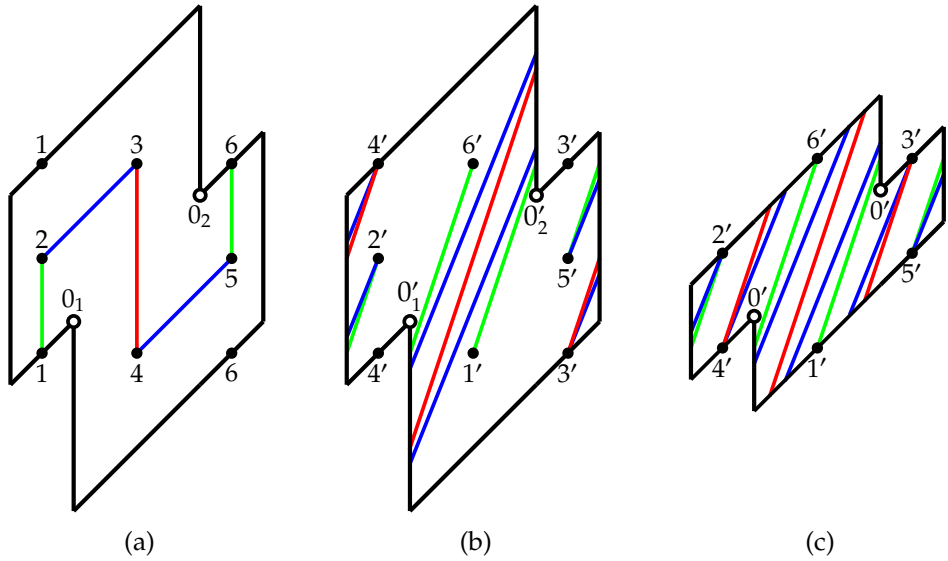
\begin{figure}
\begin{center}
  \subfigure[]{
    %\raisebox{.16\subfigwidth}{
      \resizebox{.6\subfigwidth}{!}{\begin{tikzpicture}
  \newcommand{\tw}{36pt}

  \coordinate (O) at (0,0);
  \coordinate (p1) at (-1/2*\tw,-1/2*\tw);
  \coordinate (p1b) at (-1/2*\tw,1/2*\tw);
  \coordinate (p2) at (-1/2*\tw,0);
  \coordinate (p3) at (0,1/2*\tw);
  \coordinate (p4) at (0,-1/2*\tw);
  \coordinate (p5) at (1/2*\tw,0);
  \coordinate (p6) at (1/2*\tw,1/2*\tw);
  \coordinate (p6b) at (1/2*\tw,-1/2*\tw);

  \coordinate (ll) at (-2/3*\tw,-2/3*\tw);
  \coordinate (ul) at (-2/3*\tw, 1/3*\tw);
  \coordinate (lr) at ( 1/3*\tw, 1/3*\tw);
  \coordinate (ur) at ( 1/3*\tw, 4/3*\tw);

  \begin{scope}[line width=.7pt]

    % Top torus
    \draw[thick] (ll) -- (ul);
    \draw[thick] (lr) -- (ur);
    \draw[thick] (ul) -- ($(ul) + (1/3*\tw,1/3*\tw)$);
    \draw[thick] (ll) -- ($(ll) + (1/3*\tw,1/3*\tw)$);
    \draw[thick] ($(ul) + (1/3*\tw,1/3*\tw)$) -- (ur);

    % Bottom torus
    \coordinate (sh) at (1/3*\tw,-2/3*\tw);
    \draw[thick] ($(ll)+(sh)$) -- ($(ul)+(sh)$);
    \draw[thick] ($(lr)+(sh)$) -- ($(ur)+(sh)$);
    \draw[thick] ($(ur)+(sh) - (1/3*\tw,1/3*\tw)$) -- ($(ur)+(sh)$);
    \draw[thick] ($(lr)+(sh) - (1/3*\tw,1/3*\tw)$) -- ($(lr)+(sh)$);
    \draw[thick]
      ($(ll)+(sh)$) -- ($(ll)+(sh) + (2/3*\tw,2/3*\tw)$);

    \draw[green] (p1) -- (p2);
    \draw[blue] (p2) -- (p3);
    \draw[red] (p3) -- (p4);
    \draw[blue] (p4) -- (p5);
    \draw[green] (p5) -- (p6);

    \draw[black,fill=white] ($(O) - (1/3*\tw,1/3*\tw)$) circle (.03*\tw)
      node[above,scale=.45] {\small $0_1$};
    \draw[black,fill=white] ($(O) + (1/3*\tw,1/3*\tw)$) circle (.03*\tw)
      node[below,scale=.45] {\small $0_2$};
    \filldraw[black] (p2) circle (.02*\tw)
      node[above,scale=.45] {\small $2$};
    \filldraw[black] (p1) circle (.02*\tw)
      node[below,scale=.45] {\small $1$};
    \filldraw[black] (p1b) circle (.02*\tw)
      node[above,scale=.45] {\small $1$};
    \filldraw[black] (p3) circle (.02*\tw)
      node[above,scale=.45] {\small $3$};
    \filldraw[black] (p4) circle (.02*\tw)
      node[below,scale=.45] {\small $4$};
    \filldraw[black] (p6) circle (.02*\tw)
      node[above,scale=.45] {\small $6$};
    \filldraw[black] (p6b) circle (.02*\tw)
      node[below,scale=.45] {\small $6$};
    \filldraw[black] (p5) circle (.02*\tw)
      node[below,scale=.45] {\small $5$};

  \end{scope}
\end{tikzpicture}}
    %}
    \label{fig:genus2-curve}
  }
  \hspace{1em}
  \subfigure[]{
    %\raisebox{.16\subfigwidth}{
      \resizebox{.6\subfigwidth}{!}{\begin{tikzpicture}
  \newcommand{\tw}{36pt}

  \coordinate (O) at (0,0);
  \coordinate (p1) at (0,-1/2*\tw);
  \coordinate (p2) at (-1/2*\tw,0);
  \coordinate (p3) at (1/2*\tw,1/2*\tw);
  \coordinate (p3b) at (1/2*\tw,-1/2*\tw);
  \coordinate (p4) at (-1/2*\tw,-1/2*\tw);
  \coordinate (p4b) at (-1/2*\tw,1/2*\tw);
  \coordinate (p5) at (1/2*\tw,0);
  \coordinate (p6) at (0,1/2*\tw);

  \coordinate (ppp1) at (-1/2*\tw,-3/2*\tw);
  \coordinate (ppp2) at ( 1/2*\tw, 2/2*\tw);
  \coordinate (ppp3) at (-1/2*\tw,-2/2*\tw);
  \coordinate (ppp4) at ( 1/2*\tw, 3/2*\tw);

  \coordinate (ll) at (-2/3*\tw,-2/3*\tw);
  \coordinate (ul) at (-2/3*\tw, 1/3*\tw);
  \coordinate (lr) at ( 1/3*\tw, 1/3*\tw);
  \coordinate (ur) at ( 1/3*\tw, 4/3*\tw);

  \begin{scope}[line width=.7pt]

    \coordinate (v1) at (3*\tw,9*\tw);  % green
    \coordinate (v2) at (6*\tw,15*\tw); % blue
    \coordinate (v3) at (6*\tw,18*\tw); % red

    \draw[green] (p1) -- ($(p1)+.11*(v1)$);
    \draw[green] (p2) -- ($(p2)-.055*(v1)$);
    \draw[green] (p6) -- ($(p6)-.11*(v1)$);
    \draw[green] (p5) -- ($(p5)+.055*(v1)$);

    \draw[blue] (p2) -- ($(p2)-.0275*(v2)$);
    \draw[blue] (p5) -- ($(p5)+.0275*(v2)$);
    \draw[blue] (p3b) -- ($(p3b)+.0275*(v2)$);
    \draw[blue] (p4b) -- ($(p4b)-.0275*(v2)$);
    \draw[blue] ($(ppp1)+.0275*(v2)$) -- ($(ppp2)-.0275*(v2)$);
    \draw[blue] ($(ppp3)+.0275*(v2)$) -- ($(ppp4)-.0275*(v2)$);

    \draw[red] (p4b) -- ($(p4b)-.028*(v3)$);
    \draw[red] (p3b) -- ($(p3b)+.028*(v3)$);
    \draw[red] ($-.056*(v3)$) -- ($.056*(v3)$);

    % Top torus
    \draw[thick] (ll) -- (ul);
    \draw[thick] (lr) -- (ur);
    \draw[thick] (ul) -- ($(ul) + (1/3*\tw,1/3*\tw)$);
    \draw[thick] (ll) -- ($(ll) + (1/3*\tw,1/3*\tw)$);
    \draw[thick] ($(ul) + (1/3*\tw,1/3*\tw)$) -- (ur);

    % Bottom torus
    \coordinate (sh) at (1/3*\tw,-2/3*\tw);
    \draw[thick] ($(ll)+(sh)$) -- ($(ul)+(sh)$);
    \draw[thick] ($(lr)+(sh)$) -- ($(ur)+(sh)$);
    \draw[thick] ($(ur)+(sh) - (1/3*\tw,1/3*\tw)$) -- ($(ur)+(sh)$);
    \draw[thick] ($(lr)+(sh) - (1/3*\tw,1/3*\tw)$) -- ($(lr)+(sh)$);
    \draw[thick]
      ($(ll)+(sh)$) -- ($(ll)+(sh) + (2/3*\tw,2/3*\tw)$);

    \draw[black,fill=white] ($(O) - (1/3*\tw,1/3*\tw)$) circle (.03*\tw)
      node[above,scale=.45] {$0_1'$};
    \draw[black,fill=white] ($(O) + (1/3*\tw,1/3*\tw)$) circle (.03*\tw)
      node[below,scale=.45] {$0_2'$};
    \filldraw[black] (p2) circle (.02*\tw)
      node[above,scale=.45] {\small $2'$};
    \filldraw[black] (p4) circle (.02*\tw)
      node[below,scale=.45] {\small $4'$};
    \filldraw[black] (p4b) circle (.02*\tw)
      node[above,scale=.45] {\small $4'$};
    \filldraw[black] (p6) circle (.02*\tw)
      node[above,scale=.45] {\small $6'$};
    \filldraw[black] (p1) circle (.02*\tw)
      node[below,scale=.45] {\small $1'$};
    \filldraw[black] (p3) circle (.02*\tw)
      node[above,scale=.45] {\small $3'$};
    \filldraw[black] (p3b) circle (.02*\tw)
      node[below,scale=.45] {\small $3'$};
    \filldraw[black] (p5) circle (.02*\tw)
      node[below,scale=.45] {\small $5'$};

  \end{scope}
\end{tikzpicture}}
    %}
    \label{fig:genus2-curve-1123}
  }
  \hspace{1em}
  \subfigure[]{
    \raisebox{.22\subfigwidth}{
      \resizebox{.6\subfigwidth}{!}{\begin{tikzpicture}
  \newcommand{\tw}{36pt}

  \coordinate (O) at (0,0);
  \coordinate (p1) at (0,-1/2*\tw);
  \coordinate (p2) at (-1/2*\tw,0);
  \coordinate (p3) at (1/2*\tw,1/2*\tw);
  \coordinate (p4) at (-1/2*\tw,-1/2*\tw);
  \coordinate (p5) at (1/2*\tw,0);
  \coordinate (p6) at (0,1/2*\tw);

  \coordinate (ll) at (-2/3*\tw,-2/3*\tw);
  \coordinate (ul) at (-2/3*\tw,-1/6*\tw);
  \coordinate (lr) at ( 1/3*\tw, 1/3*\tw);
  \coordinate (ur) at ( 1/3*\tw, 5/6*\tw);

  \begin{scope}[line width=.7pt]

%    \draw[densely dotted]
%      ($(ll) + (1/3*\tw,1/3*\tw)$) -- ($(ul) + (1/3*\tw,1/3*\tw)$);
%    \coordinate (sh) at (1/3*\tw,-1/6*\tw);
%    \draw[densely dotted]
%      ($(lr)+(sh) - (1/3*\tw,1/3*\tw)$) -- ($(ur)+(sh) - (1/3*\tw,1/3*\tw)$);

    \coordinate (v1) at (3*\tw,9*\tw);  % green
    \coordinate (v2) at (6*\tw,15*\tw); % blue
    \coordinate (v3) at (6*\tw,18*\tw); % red

    \draw[green] (p1) -- ($(p1)+.11*(v1)$);
    \draw[green] (p2) -- ($(p2)-.055*(v1)$);
    \draw[green] (p6) -- ($(p6)-.11*(v1)$);
    \draw[green] (p5) -- ($(p5)+.055*(v1)$);

    \draw[blue] (p2) -- ($(p2)-.0275*(v2)$);
    \draw[blue] (p5) -- ($(p5)+.0275*(v2)$);
    \draw[blue] (p3) -- ($(p3)-.055*(v2)$);
    \draw[blue] (p4) -- ($(p4)+.055*(v2)$);
    \draw[blue] ($(ppp1)+.055*(v2)$) -- ($(ppp2)-.0275*(v2)$);
    \draw[blue] ($(ppp3)+.0275*(v2)$) -- ($(ppp4)-.055*(v2)$);

    \draw[red] (p4) -- ($(p4)+.042*(v3)$);
    \draw[red] (p3) -- ($(p3)-.042*(v3)$);
    \draw[red] ($-.042*(v3)$) -- ($.042*(v3)$);

    % Top torus
    % Periodic edges
    \draw[thick] (ll) -- (ul);
    \draw[thick] (lr) -- (ur);
    % Puncture 1
    \draw[thick] (ul) -- ($(ul) + 1/2*(1/3*\tw,1/3*\tw)$);
    \draw[thick]
      ($(ul) + (1/3*\tw,1/3*\tw)$) -- ($(ul) + 1/2*(1/3*\tw,1/3*\tw)$);
    % Puncture 2
    \draw (ll) -- ($(ll) + 1/2*(1/3*\tw,1/3*\tw)$);
    \draw
      ($(ll) + (1/3*\tw,1/3*\tw)$) -- ($(ll) + 1/2*(1/3*\tw,1/3*\tw)$);
    % Puncture 3
    \draw
      ($(ul) + (1/3*\tw,1/3*\tw)$) -- ($(ul) + (2/3*\tw,2/3*\tw)$);
    \draw
      (ur) -- ($(ul) + (2/3*\tw,2/3*\tw)$);

    % Bottom torus
    \coordinate (sh) at (1/3*\tw,-1/6*\tw);
    % Periodic edges
    \draw[thick] ($(ll)+(sh)$) -- ($(ul)+(sh)$);
    \draw[thick] ($(lr)+(sh)$) -- ($(ur)+(sh)$);
    % Puncture 5
    \draw[thick] ($(ur)+(sh)$) -- ($(ur)+(sh) - 1/2*(1/3*\tw,1/3*\tw)$);
    \draw[thick]
      ($(ur)+(sh) - (1/3*\tw,1/3*\tw)$) -- ($(ur)+(sh) - 1/2*(1/3*\tw,1/3*\tw)$);   % Puncture 6
    \draw[thick] ($(lr)+(sh)$) -- ($(lr)+(sh) - 1/2*(1/3*\tw,1/3*\tw)$);
    \draw[thick]
      ($(lr)+(sh) - (1/3*\tw,1/3*\tw)$) -- ($(lr)+(sh) - 1/2*(1/3*\tw,1/3*\tw)$);
    % Puncture 4
    \draw
      ($(ll)+(sh)$) -- ($(ll)+(sh) + (1/3*\tw,1/3*\tw)$);
    \draw
      ($(ll)+(sh) + (2/3*\tw,2/3*\tw)$) -- ($(ll)+(sh) + (1/3*\tw,1/3*\tw)$);

    \draw[black,fill=white] ($(O) - (1/3*\tw,1/3*\tw)$) circle (.03*\tw)
      node[above,scale=.45] {\small $\,\,0'$};
    \draw[black,fill=white] ($(O) + (1/3*\tw,1/3*\tw)$) circle (.03*\tw)
      node[below,scale=.45] {\small $0'$};

    \filldraw[black] (p1) circle (.02*\tw)
      node[below,scale=.45] {\small $1'$};
    \filldraw[black] (p2) circle (.02*\tw)
      node[above,scale=.45] {\small $2'$};
    \filldraw[black] (p3) circle (.02*\tw)
      node[above,scale=.45] {\small $3'$};
    \filldraw[black] (p4) circle (.02*\tw)
      node[below,scale=.45] {\small $4'$};
    \filldraw[black] (p5) circle (.02*\tw)
      node[below,scale=.45] {\small $5'$};
    \filldraw[black] (p6) circle (.02*\tw)
      node[above,scale=.45] {\small $6'$};

  \end{scope}
\end{tikzpicture}}
    }
    \label{fig:genus2_sphere-curve-1123}
  }
\end{center}
\caption{(a) The genus two surface from Fig.~\ref{fig:genus2_flat}, with
  opposite edges identified and arcs between the rods.  (b) The surface and
  arcs after applying the map~\eqref{eq:1123} and using the edge
  identifications to cut up and rearrange the surface to the same initial
  domain.  (c) The arcs on the punctured sphere~$\S_{0,6}=\S_2/\iota_2$, with
  edges identified as in Fig.~\ref{fig:genus2_sphere}.}
\label{fig:genus2-1123}
\end{figure}

In Fig.~\ref{fig:genus2-curve} we reproduce the genus two surface, omitting
the edge identifications for clarity, and draw some arcs between our rods.
Now act on the surface (embedded in the plane) with the map~\eqref{eq:1123}.
The polygon gets stretched, and we cut and glue pieces following the edge
identifications to bring it back into its initial domain, as in
Fig.~\ref{fig:genus2-curve-1123}.  Punctures~$2$ and~$5$ are fixed, $1$ and
$4$ are swapped, as are $3$ and $6$.  This is exactly the same as for a
half-period of the puller in Fig.~\ref{fig:6rods_rodmotion}.

\begin{figure}
\begin{center}
  \subfigure[]{
    \begin{tikzpicture}
  \newcommand{\tw}{36pt}
  \coordinate (q0) at (-\tw,0);
  \coordinate (q1) at (0,0);
  \coordinate (q2) at (\tw,0);
  \coordinate (q3) at (2*\tw,0);
  \coordinate (O) at (2*\tw,0);
  \coordinate (q4) at (3*\tw,0);
  \coordinate (q5) at (4*\tw,0);
  \coordinate (q6) at (5*\tw,0);

  \begin{scope}[line width=1.2pt]
    % Draw outer ellipse, to suggest a sphere.
    \tikzbgellipsesix

    \draw[green,ultra thick]%,postaction={on each segment={mid arrow}}]
      (q1) to (q2);

    \draw[blue,ultra thick]%,postaction={on each segment={mid arrow}}]
      (q2) to (q3);

    \draw[red,ultra thick]%,postaction={on each segment={mid arrow}}]
      (q3) to (q4);

    \draw[blue,ultra thick]%,postaction={on each segment={mid arrow}}]
      (q4) to (q5);

    \draw[green,ultra thick]%,postaction={on each segment={mid arrow}}]
      (q5) to (q6);

    \filldraw[black,fill=white] (q0) circle (.03*\tw)
      node[below] {\small $0$};
    \filldraw[black] (q1) circle (.03*\tw)
      node[below] {\small $1$};
    \filldraw[black] (q2) circle (.03*\tw)
      node[below] {\small $2$};
    \filldraw[black] (q3) circle (.03*\tw)
      node[below] {\small $3$};
    \filldraw[black] (q4) circle (.03*\tw)
      node[below] {\small $4$};
    \filldraw[black] (q5) circle (.03*\tw)
      node[below] {\small $5$};
    \filldraw[black] (q6) circle (.03*\tw)
      node[below] {\small $6$};
  \end{scope}
\end{tikzpicture}
  }
  \subfigure[]{
    \begin{tikzpicture}
  \newcommand{\tw}{36pt}
  \coordinate (q0) at (-\tw,0);
  \coordinate (q4) at (0,0);
  \coordinate (q2) at (\tw,0);
  \coordinate (q6) at (2*\tw,0);
  \coordinate (O) at (2*\tw,0);
  \coordinate (q1) at (3*\tw,0);
  \coordinate (q5) at (4*\tw,0);
  \coordinate (q3) at (5*\tw,0);

  \begin{scope}[line width=1.2pt]
    % Draw outer ellipse, to suggest a sphere.
    \tikzbgellipsesix

    \draw[green,ultra thick]%,postaction={on each segment={mid arrow}}]
      (q1)
      to[out=110,in=70] (q2)
      ;

    \draw[blue,ultra thick]%,postaction={on each segment={mid arrow}}]
      (q2)
      to[out=40,in=130] ($(q1) - (.25*\tw,0)$)
      to[out=-50,in=-140] ($(q1) + (.25*\tw,0)$)
      to[out=40,in=130] (q3)
      ;

    \draw[red,ultra thick]%,postaction={on each segment={mid arrow}}]
      (q3)
      to[out=150,in=35] ($(q5) - (.5*\tw,0)$)
      to[out=-145,in=-45] ($(q6) + (.5*\tw,0)$)
      to[out=135,in=35] ($(q2) + (.5*\tw,0)$)
      to[out=-145,in=-30] (q4)
      ;

    \draw[green,ultra thick]%,postaction={on each segment={mid arrow}}]
      (q6)
      to[out=-70,in=-110] (q5)
      ;

    \draw[blue,ultra thick]%,postaction={on each segment={mid arrow}}]
      (q5)
      to[out=-130,in=-40] ($(q6) + (.25*\tw,0)$)
      to[out=140,in=50] ($(q6) - (.25*\tw,0)$)
      to[out=-140,in=-50] (q4)
      ;

    \filldraw[black,fill=white] (q0) circle (.03*\tw)
      node[below] {\small $0'$};
    \filldraw[black] (q1) circle (.03*\tw)
      node[above right] {\small $1'$};
    \filldraw[black] (q2) circle (.03*\tw)
      node[above left] {\small $2'$};
    \filldraw[black] (q3) circle (.03*\tw)
      node[below] {\small $3'$};
    \filldraw[black] (q4) circle (.03*\tw)
      node[below left] {\small $4'$};
    \filldraw[black] (q5) circle (.03*\tw)
      node[below right] {\small $5'$};
    \filldraw[black] (q6) circle (.03*\tw)
      node[below left] {\small $6'$};
  \end{scope}
\end{tikzpicture}
  }
\end{center}
\caption{(a) A sphere with six punctures (rods) and a seventh puncture at the
  fixed point~$0$, with arcs between the punctures, as in
  Fig.~\ref{fig:genus2-curve}.  (b) The arcs after a half-period of the rod
  motion in Fig.~\ref{fig:6rods_rodmotion}.  These are identical to the arcs
  of Fig~\ref{fig:genus2_sphere-curve-1123}.}
\label{fig:genus2_sphere-curve}
\end{figure}
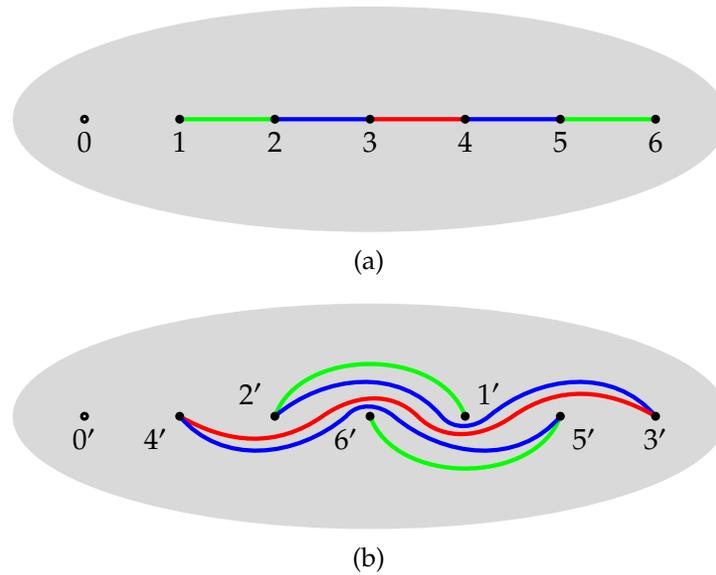

After acting with the map we form the quotient
surface~$\S_2/\iota_2=\S_{0,6}$, as in
Fig.~\ref{fig:genus2_sphere-curve-1123}.  Now we can carefully trace out the
path of each arc, and keep track of which side of the arcs the punctures lie.
The paths in Fig.~\ref{fig:genus2_sphere-curve-1123} are identical to the arcs
in Fig.~\ref{fig:genus2_sphere-curve}, and we conclude that the
map~\eqref{eq:1123} is the correct description of the six-rod puller.  Its
dilatation is thus the largest root of~$x^2 - 4x +1$, which is~$2+\sqrt{3}$.

The description of the surface as a polygon in the plane, with edge
identifications via translations and rotations, comes from the theory of
\emph{flat surfaces} \autocite{Zorich2006}.  In this viewpoint the surface is
given a flat metric, and the corners of the polygon correspond to
\emph{conical singularities} with infinite curvature.  Here, the two
singularities~$0_{1,2}$ have cone angle~$4\pi$, as can be seen by drawing a
small circle around the points and following the edge identifications.  The
sum of the two singularities is~$8\pi$, which equals~$2\pi(4g-4)$ by the
Gauss--Bonnet formula, with~$g=2$ the genus.

\keepnote{Discuss commutation of~$\phi$ with $\iota_{1,2}$.  What is~$-I$?  Is
  it~$\iota_1\iota_2$?}

\section*{What is the best taffy puller?}

There are many other taffy puller designs found in the patent literature.%
\ifarxiv{ (See the Appendix for some examples.)}{} A few of these have a
quadratic dilation, like the examples we discussed, but many don't: they
involve pseudo-Anosov maps that are more complicated than simple branched
covers of the torus.  We will not give a detailed construction of the maps,
but rather report the polynomial whose largest root is the dilatation and
offer some comments.  The polynomials were obtained using the computer
programs \texttt{braidlab} \autocite{braidlab} and \texttt{train}
\autocite{HallTrain}.

Many taffy pullers are planetary devices --- these have rods that move on
epicycles, giving their orbits a `spirograph' appearance.  The name comes from
Ptolemaic models of the solar system, where planetary motions were apparently
well-reproduced using systems of gears.  Planetary designs are used in many
mixing devices, and are a natural way of creating taffy pullers.
\textcite{Kobayashi2007,Kobayashi2010} and \textcite{MattFinn2011_silver} have
designed and studied a class of such devices.

\begin{figure}
\vspace{.25in}
\begin{center}
\subfigure[]{
  \includegraphics[width=.29\textwidth]{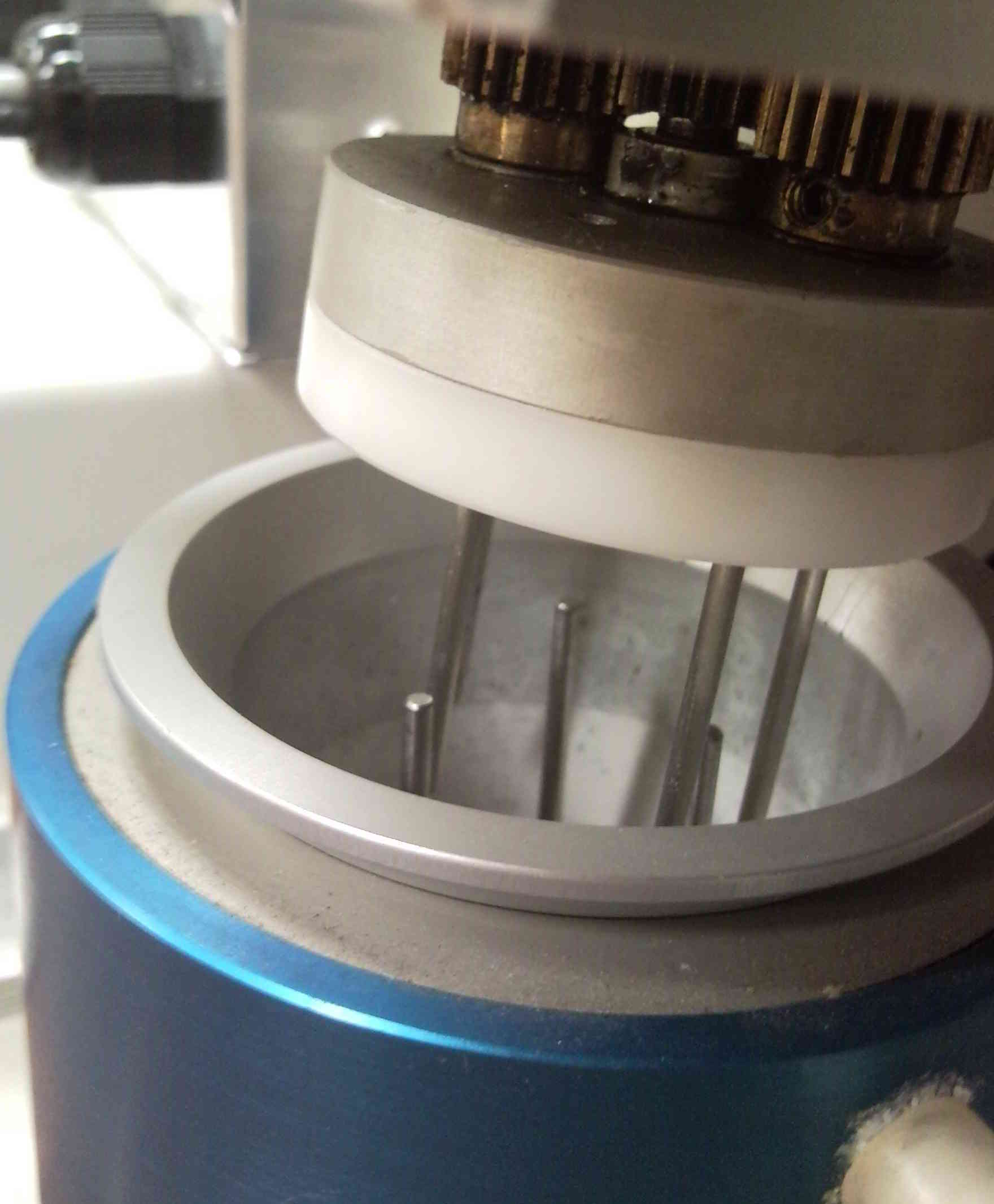}
  \label{fig:mixograph6}
}\hspace{1em}%
\subfigure[]{
  \begin{minipage}{.183\textwidth}
    \vspace{-1.8in}
    \includegraphics[width=\textwidth]{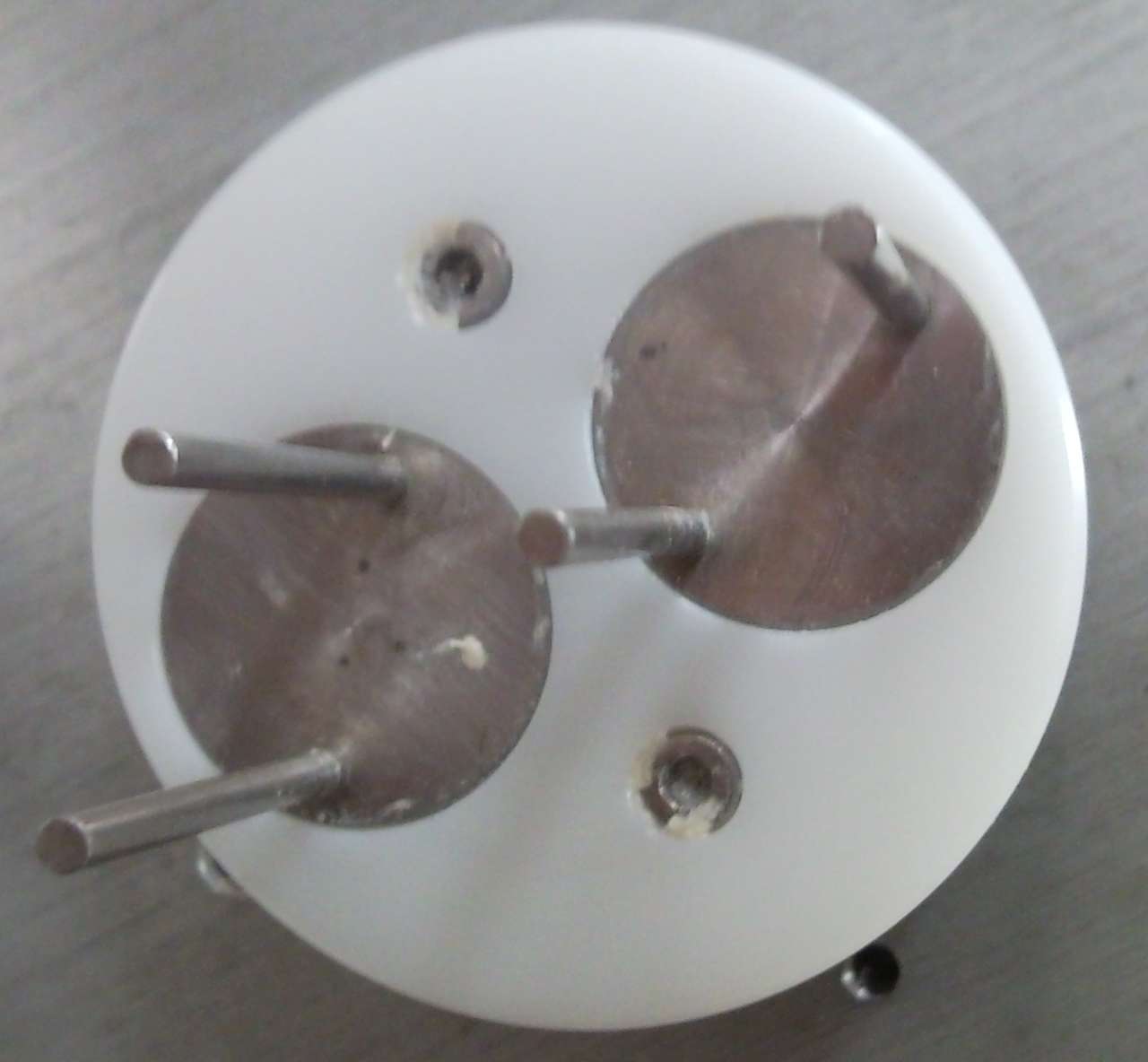}
    \includegraphics[width=\textwidth]{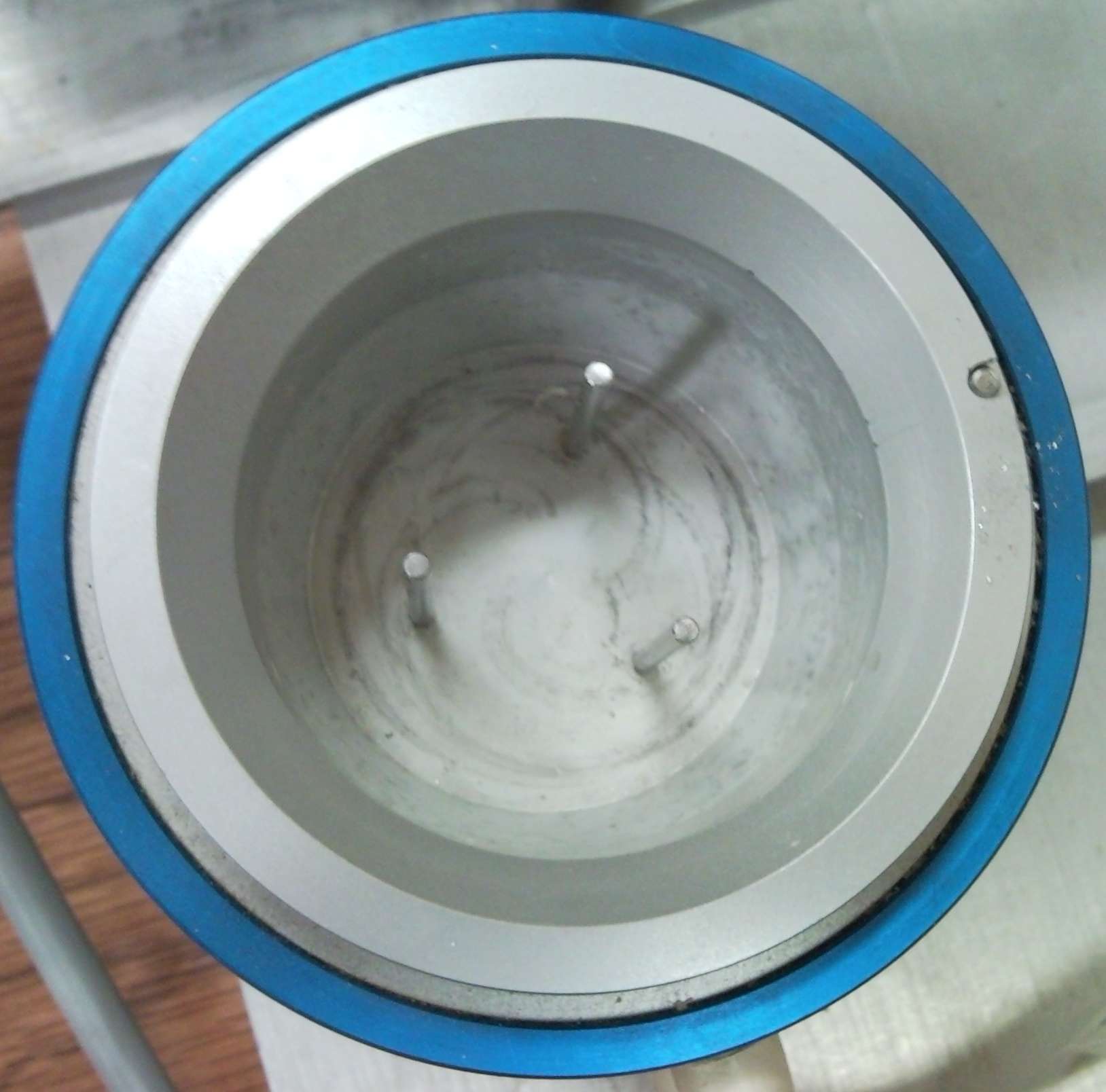}
    \label{fig:mixograph34}
  \end{minipage}
}\hspace{1em}%
\subfigure[]{
  \begin{minipage}{.25\textwidth}
    \vspace{-1.8in}
    \includegraphics[width=\textwidth]{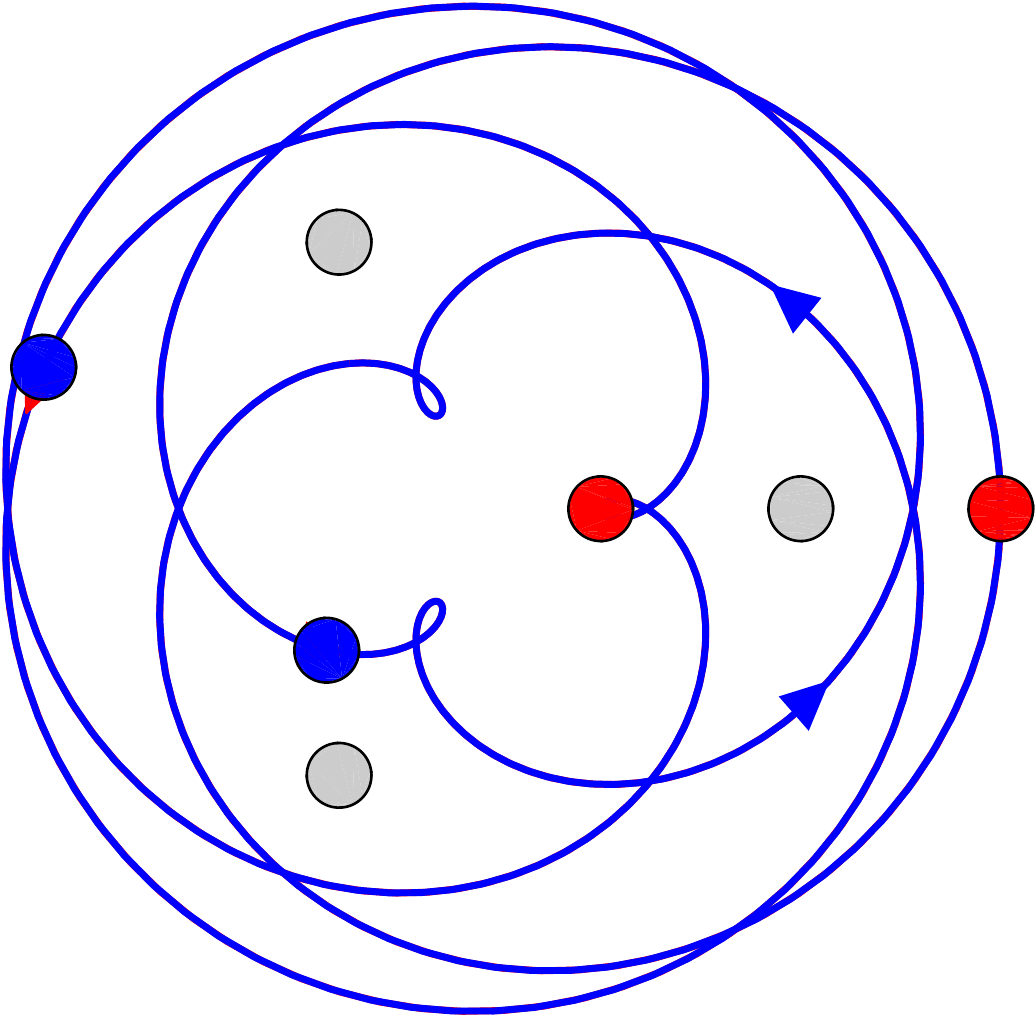}

    \vspace{.5em}

    \includegraphics[width=\textwidth]{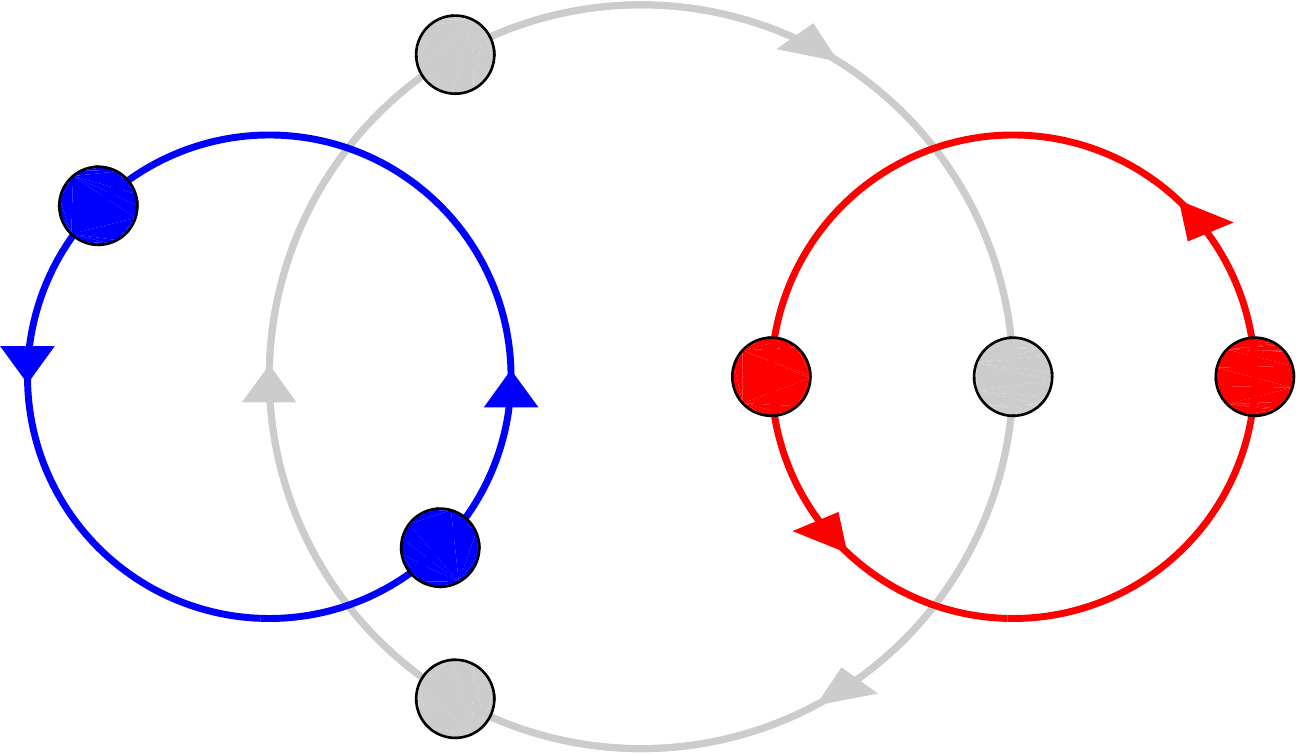}
    \label{fig:mixograph_rodmotion}
  \end{minipage}
}
\end{center}
\caption{(a) The mixograph, a planetary rod mixer for bread dough. (b) Top
  section with four moving rods (above), and bottom section with three fixed
  rods (below). (c) The rod motion is complex (top), but is less so in a
  rotating frame (bottom).  (Courtesy of the Department of Food Science,
  University of Wisconsin.  Photo by the author, from
  \textcite{MattFinn2011_silver}.)}
\label{fig:mixograph}
\end{figure}

A typical planetary device, the mixograph, is shown in
Fig.~\ref{fig:mixograph}.  The mixograph consists of a small cylindrical
vessel with three fixed vertical rods.  A lid is lowered onto the base.  The
lid has two gears each with a pair of rods, and is itself rotating, resulting
in a net complex motion as in Fig.~\ref{fig:mixograph_rodmotion} (top).  The
mixograph is used to measure properties of bread dough: a piece of dough is
placed in the device, and the torque on the rods is recorded on graph paper,
in a similar manner to a seismograph.  An expert on bread dough can then
deduce dough-mixing characteristics from the graph \autocite{Connelly2008}.

Clearly, passing to a uniformly-rotating frame does not modify the dilatation.
For the mixograph, a co-rotating frame where the fixed rods rotate simplifies
the orbits somewhat (Fig.~\ref{fig:mixograph_rodmotion}, bottom).
\keepnote{This has another advantage: in this `co-rotating' picture, the rods
  return to their initial configuration (as a set) much sooner than in the
  fixed picture.  This makes the map simpler, and we then take a power of this
  co-rotating map to recover the original.}  The rod motion of
Fig.~\ref{fig:mixograph_rodmotion} (bottom) must be repeated \emph{six} times
for all the rods to return to their initial position.  The dilatation for the
co-rotating map is the largest root of
$x^8 - 4x^7 - x^6 + 4x^4 - x^2 - 4x + 1$, which is approximately $4.1858$.

%\begin{landscape}
\begin{table}
  \caption{Efficiency of taffy pullers.  A number of rods such as~6 (2)
    indicates~$6$ total rods, with~$2$ fixed.  The largest root of the
    polynomial is the dilatation.  The dilatation corresponds to a fraction~$p$
    of a full period, when each rod returns to its initial position.
    The entropy per period is
    $\log(\text{dilatation})/\text{period}$, which is a crude measure of
    efficiency.  Here $\varphi = \tfrac12(1+\sqrt5)$ is the Golden Ratio, and
    $\chi = 1 + \sqrt2$ is the Silver Ratio.}
\label{tab:taffypull}
{\small
\begin{center}
\begin{tabular}{lllllll}
\hline\hline
puller & fig. & rods & polynomial & dilatation & $p$ &
\begin{tabular}{@{}c@{}}entropy/\\ period\end{tabular} \\[2pt]
\hline
standard 3-rod & \ref{fig:Robinson1908} & 3 (1) & $x^2-6x+1$
  & $\chi^2$ & 1 & $1.7627$ \\
Nitz (1918) & \ref{fig:s1s-2} & 3 & $x^2-3x+1$ & $\varphi^2$ &
  \nicefrac{1}{3} & $2.8873$ \\
standard 4-rod & \ref{fig:Thibodeau1903} & 4 &
  $x^2-6x+1$ & $\chi^2$ & 1 & $1.7627$ \\
Thibodeau (1904) & \ifarxiv{\ref{fig:Thibodeau1904}}{--} & 4 & $x^2-3x+1$ &
  $\varphi^2$ & \nicefrac{1}{3} & $2.8873$ \\
6-rod & \ref{fig:6rods} & 6 (2) & $x^2-4x+1$ & $2+\sqrt3$ & \nicefrac{1}{2} &
  $2.6339$ \\
McCarthy (1916)\,\dag &
  \ifarxiv{\ref{fig:McCarthy1916a_rodmotion}}{--} & 4 (3)
  & $x^2-18x+1$ & $\varphi^6$ & 1 & $2.8873$ \\
  & \ifarxiv{\ref{fig:McCarthy1916b_rodmotion}}{--} & 4 (3) &
  $x^4-36x^3+54x^2-36x+1$ & $34.4634$ & 1 & $3.5399$ \\
mixograph\,\ddag & \ref{fig:mixograph_rodmotion} & 7 &
  \begin{tabular}{@{}l@{}}$x^8 - 4x^7 - x^6 + 4x^4$\\$\ \ \ \, \qquad\qquad\hbox{} - x^2 - 4x + 1$\end{tabular} & $4.1858$ & \nicefrac{1}{6}
  & $8.5902$ \\
%mixograph\,\ddag & \ref{fig:mixograph_rodmotion} & 7 &
%  $x^8 - 4x^7 - x^6 + 4x^4 - x^2 - 4x + 1$ & $4.1858$ & \nicefrac{1}{6}
%  & $8.5902$ \\
Jenner (1905) & \ifarxiv{\ref{fig:Jenner1905}}{--}
  & 5 (3) & $x^4-8x^3-2x^2-8x+1$ &
  $(\varphi+\sqrt{\varphi})^2$ & 1 & $2.1226$ \\
Shean (1914) & \ifarxiv{\ref{fig:Shean1914}}{--}
  & 6 & $x^2-4x+1$ & $2+\sqrt3$ &
  \nicefrac{1}{2} & $2.6339$ \\
McCarthy (1915) & \ifarxiv{\ref{fig:McCarthy1915}}{--} & 5 (2) &
  $x^4 - 20x^3 - 26x^2 - 20x + 1$ & $21.2667$ & 1 &
  $3.0571$ \\
\hline\hline
\end{tabular}
\end{center}
\begin{flushleft}
\dag\, The McCarthy (1916) device has two configurations. \\
\ddag\, This is the co-rotating version of the mixograph
  (Fig.~\ref{fig:mixograph_rodmotion}, bottom).
\end{flushleft}
}
\end{table}
%\end{landscape}

The reader might be wondering at this point: which is the best taffy puller?
Did all these incremental changes and new designs in the patent literature
lead to measurable progress in the effectiveness of taffy pullers?
Table~\ref{tab:taffypull} collects the characteristic polynomials and the
dilatations (the largest root) for all the taffy pullers discussed here%
\ifarxiv{ and a few others included at the end}{ and a few others}%
.  The total number of rods is listed (the number in parentheses is the number
of fixed rods).

The column labeled~$p$ requires a bit of explanation.  Comparing the different
taffy pullers is not straightforward.  To keep things simple, we take the
efficiency to be the total dilatation for a full period, defined by all the
rods returning to their initial position.  For example, referring to
Table~\ref{tab:taffypull}, for the \cite{nitz_candy-puller._1918} device the
rods return to the same configuration (as a set) after~$p=1/3$ period.  Hence,
the dilatation listed, $\varphi^2$, is for~$1/3$ period.  We define a puller's
efficiency as the entropy (logarithm of the dilatation) per period.  In this
case the efficiency is $\log(\varphi^2)/(1/3) = 6\log\varphi \approx 2.8873$.
By this measure, the mixograph is the clear winner, with a staggering
efficiency of~$8.5902$.  Of course, it also has the most rods.  The large
efficiency is mostly due to how long the rods take to return to their initial
position.

Some general observations can be made regarding practical taffy pullers.  With
a few exceptions, they all give pseudo-Anosov maps.  Though we did not define
this term precisely, in this context it implies that any initial piece of
taffy caught on the rods will grow exponentially.  The inventors were thus
aware, at least intuitively, that there should be no unnecessary rods.
Another observation is that most of the dilatations are quadratic numbers.
There are probably a few reasons for this.  One is that the polynomial giving
the dilatation expresses a recurrence relation that characterizes how the
taffy's folds are combined at each period.  With a small number of rods, there
is a limit to the degree of this recurrence ($2n-4$ for $n$ rods).  A second
reason is that more rods does not necessarily mean larger dilatation
\autocite{MattFinn2011_silver}.  On the contrary, more rods allows for a
smaller dilatation, as observed when finding the smallest value of the
dilatation \autocite{Hironaka2006, Thiffeault2006, Venzke_thesis,
  LanneauThiffeault2011_braids}.

The collection of taffy pullers presented here can be thought of as a battery
of examples to illustrate various types of pseudo-Anosov maps.  Even though
they did not come out of the mathematical literature, they predate by many
decades the examples that were later constructed by mathematicians
\autocite{Boyland2000, Thiffeault2006, Kobayashi2007, Binder2008, Binder2010,
  Kobayashi2010, MattFinn2011_silver, Boyland2011}.

\subsection*{Acknowledgments} The author thanks Alex Flanagan for
helping to design and build the 6-rod taffy puller, and Phil Boyland and Eiko
Kin for their comments on the manuscript.

\printbibliography

% Do not remove the tags BEGIN_ARXIV / END_ARXIV.

%% BEGIN_ARXIV

\ifarxiv{}{\end{document}}

%\newpage

\section*{Appendix: A few more taffy puller designs}

Because of space constraints, several taffy pullers from the patent literature
were omitted from the
\href{http://dx.doi.org/10.1007/s00283-018-9788-4}{\emph{Mathematical
    Intelligencer}} version of this article.  We include these here for the
interested reader.

\subsection*{A 4-rod device with Golden ratio dilatation.}
\begin{figure}
\begin{center}
\begin{minipage}{.5\textwidth}
\begin{center}
\subfigure[]{
  \includegraphics[width=\textwidth]{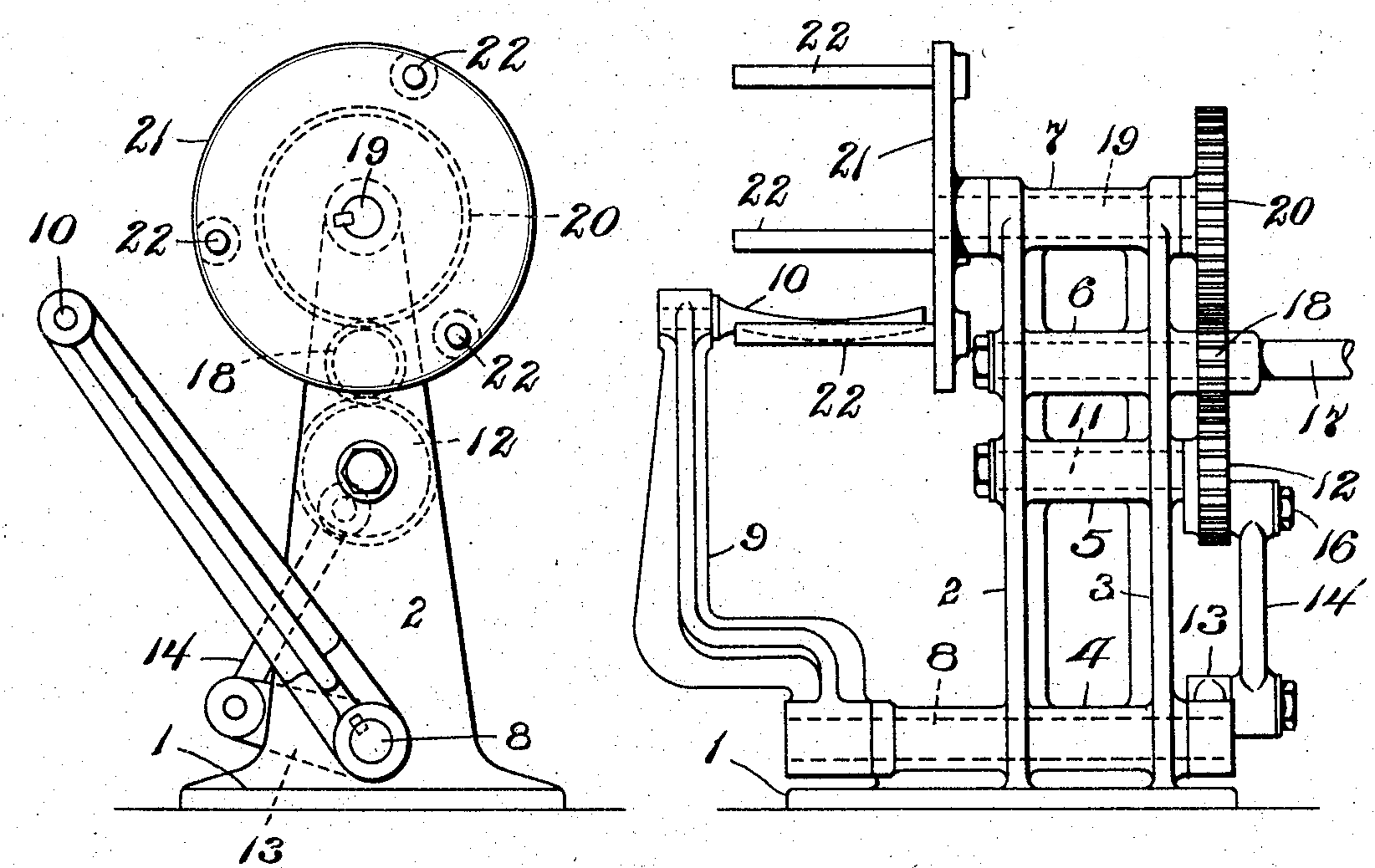}
  \label{fig:Thibodeau1904_device}
}\hspace{1em}
\subfigure[]{
  \raisebox{0em}{
  \includegraphics[width=.55\textwidth]{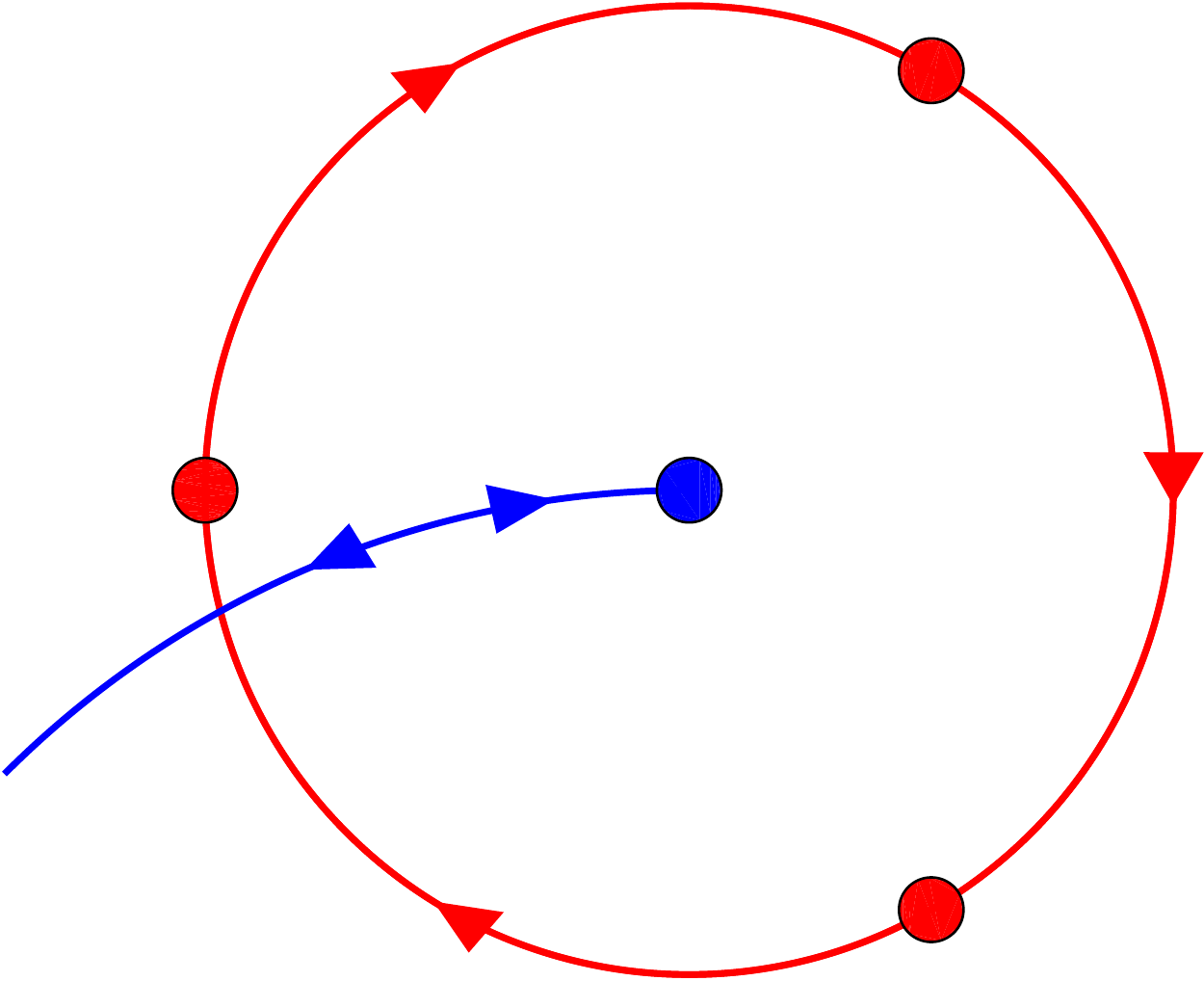}}
  \label{fig:Thibodeau1904_rodmotion}
}
\end{center}
\end{minipage}
\hspace{1em}
\begin{minipage}{.38\textwidth}
\subfigure[]{
  \includegraphics[width=\textwidth]{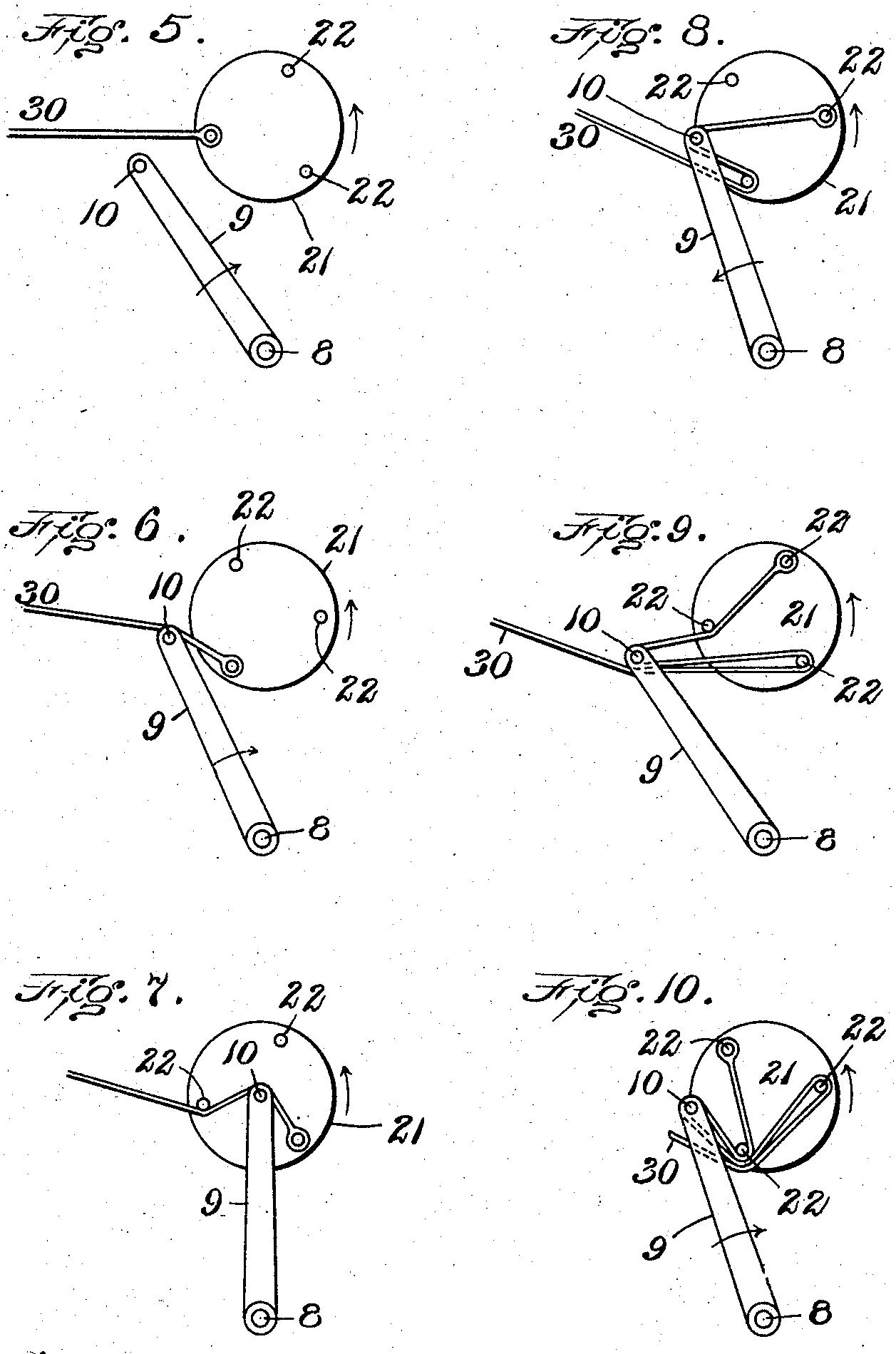}
  \label{fig:Thibodeau1904_action}
}
\end{minipage}
\end{center}
\caption{(a) Taffy puller from the patent of
  \textcite{thibodeau_candy-pulling_1904}, with three rotating rods on a wheel
  and an oscillating arm.  (b) Rod motion.  (c) The action of the taffy
  puller, as depicted in the patent.}
\label{fig:Thibodeau1904}
\end{figure}

Thibodeau's device in Fig.~\ref{fig:Thibodeau1903} gave an example of a 4-rod
taffy puller arising directly from an Anosov map.  Another example is a later
design of Thibodeau shown in Fig.~\ref{fig:Thibodeau1904}.  It consists of
three rods moving in a circle, and a fourth rod crossing their path back and
forth.  This taffy puller can be shown to come from the same Anosov as gave us
the 3-rod puller in Fig.~\ref{fig:s1s-2}, with dilatation equal to the square
of the Golden Ratio.  Thus, if one is interested in building a device with a
Golden Ratio dilatation, the design in Fig.~\ref{fig:Thibodeau1904_device} is
probably far easier to implement than Nitz's in
Fig.~\ref{fig:Nitz1918_device}, since Thibodeau's does not involve rods being
exchanged between two gears.

\subsection*{A simple planetary design.}
\begin{figure}
\begin{center}
\subfigure[]{
  \includegraphics[height=.3\textheight]{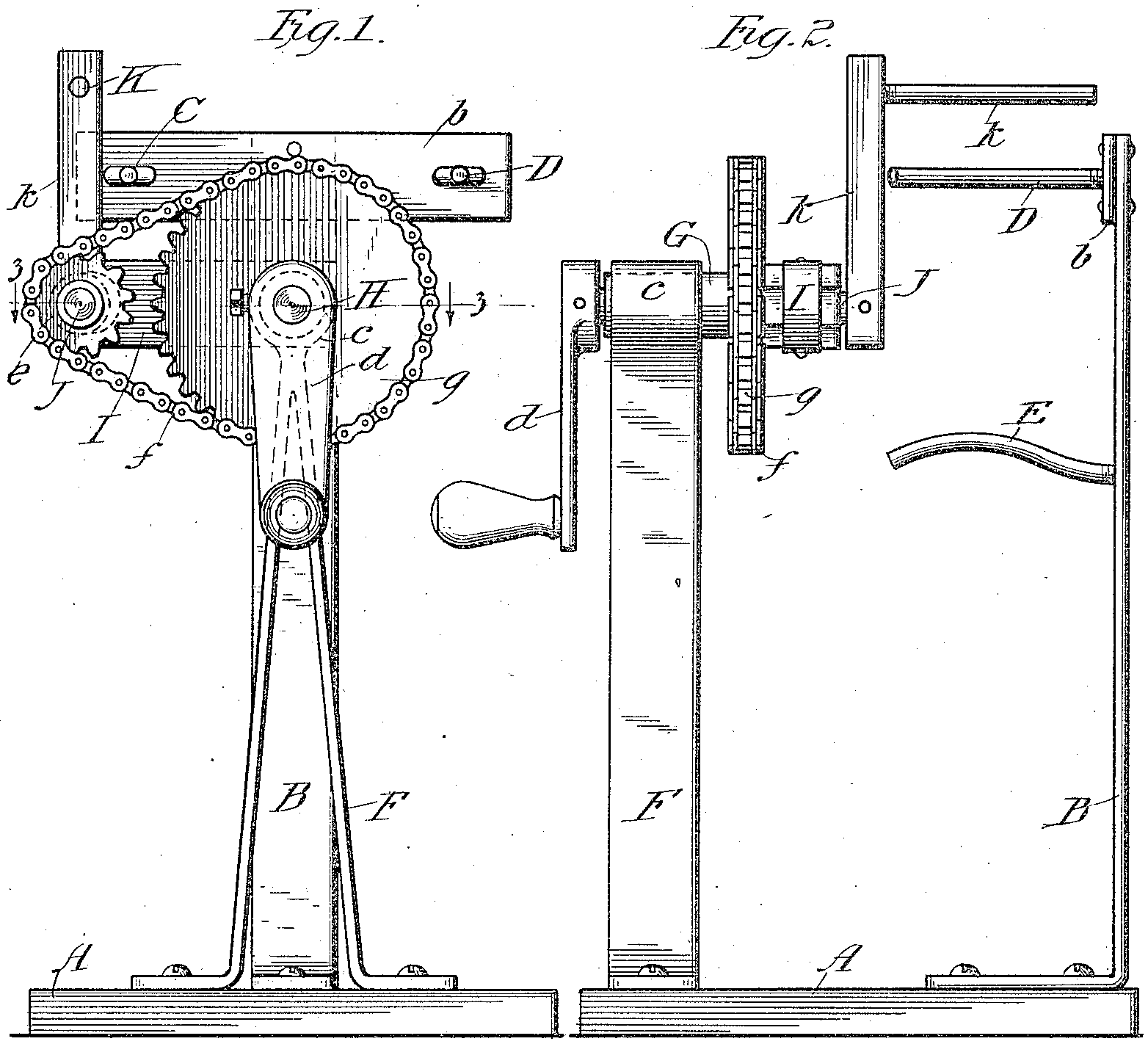}
  \label{fig:McCarthy1916_device}
}\hspace{1em}
\subfigure[]{
  \includegraphics[height=.3\textheight]{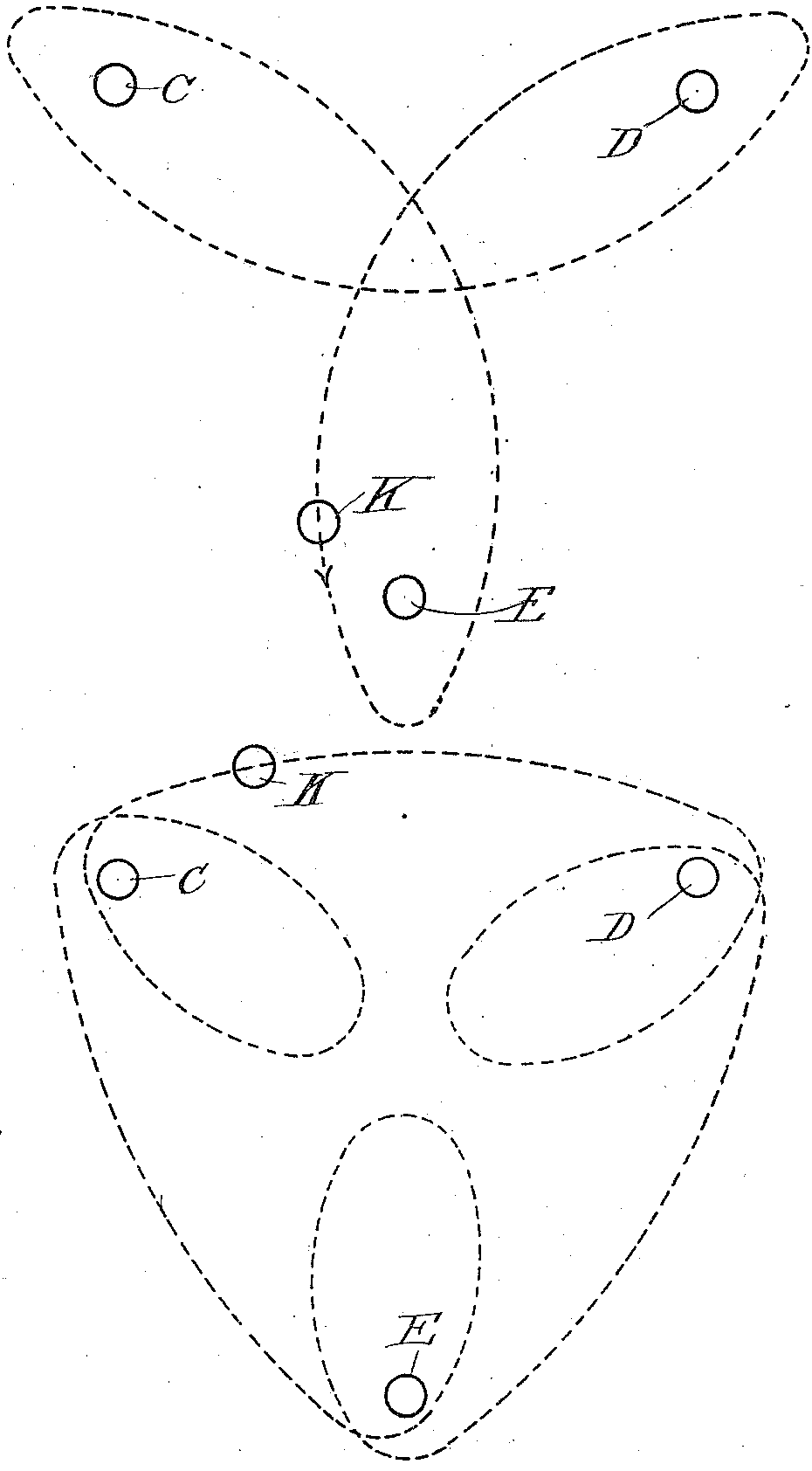}
  \label{fig:McCarthy1916_rodmotion2}
}

\subfigure[]{
  \includegraphics[height=.2\textheight]{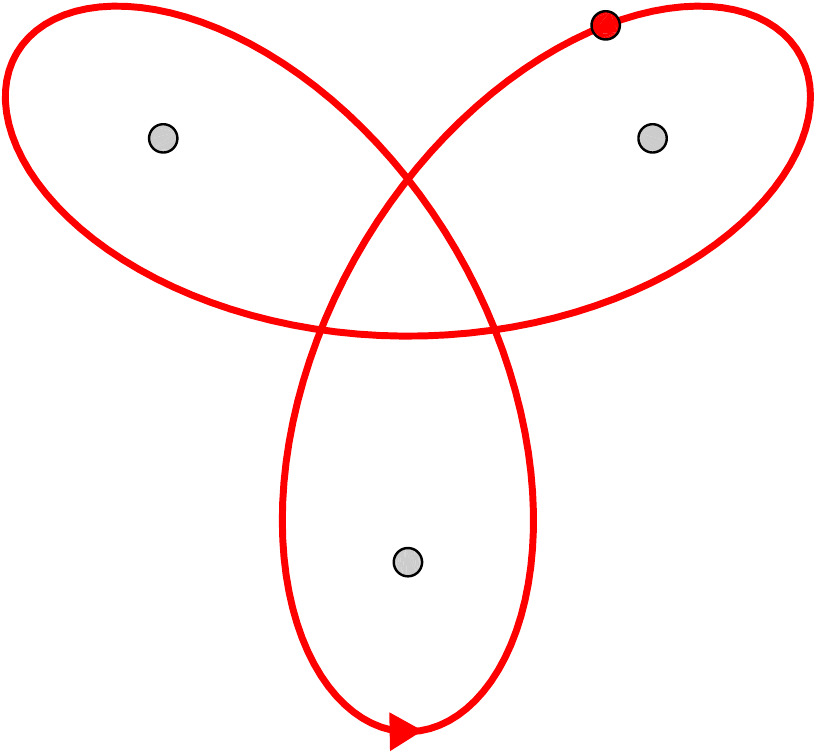}
  \label{fig:McCarthy1916a_rodmotion}
}\hspace{1em}
\subfigure[]{
  \includegraphics[height=.2\textheight]{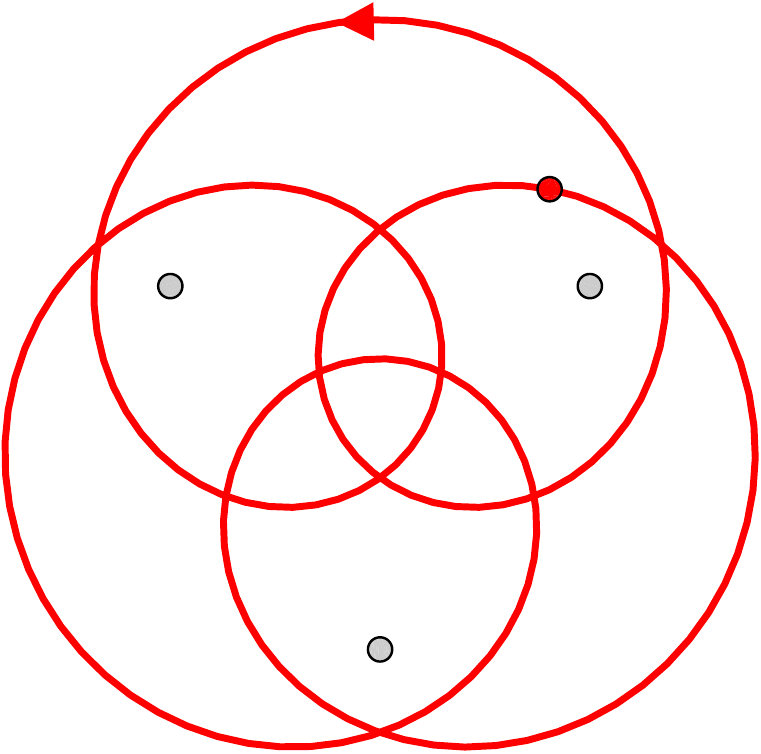}
  \label{fig:McCarthy1916b_rodmotion}
}
\end{center}
\caption{(a) Taffy puller from the patent of
  \textcite{mccarthy_candy-pulling_1916}.  (b)~Rod motion for the two
  configurations of the device as sketched in the patent.  (c) and (d): the
  actual rod motions.}
\label{fig:McCarthy1916}
\end{figure}

\textcite{mccarthy_candy-pulling_1916} has an interesting planetary design for
a taffy puller (Fig.~\ref{fig:McCarthy1916}).  It has two configurations, with
rod motions shown in Fig.~\ref{fig:McCarthy1916_rodmotion2}.  Its first
configuration (pictured in Fig.~\ref{fig:McCarthy1916_device} with rod motion
as in Fig.~\ref{fig:McCarthy1916a_rodmotion}) is a perfect example of a
`$\pi_1$-stirring device,' a device where only a single rod moves around a set
of fixed rods.  The optimality of such devices was studied by
\textcite{Boyland2011}, and McCarthy's device is one of their optimal
examples.

The second configuration (not shown) involves replacing the chain in
Fig.~\ref{fig:McCarthy1916_device} by two gears in direct contact.  This gives
the motion in Fig.~\ref{fig:McCarthy1916b_rodmotion}, which does appear quite
different from McCarthy's sketch (Fig.~\ref{fig:McCarthy1916_rodmotion2},
bottom) but is topologically identical.  McCarthy himself seemed to prefer the
first configuration, as he noted a bit wordily in his patent:
\begin{quote}
  The planetary course described by this pin, when this modified construction
  is employed, gives a constant pull to the candy, but does not accomplish as
  thorough mixing of the same as when said pin describes the planetary course
  resulting from the construction of the preferred form of my invention, as
  hereinbefore first described.
\end{quote}
What he meant by `a constant pull to the candy' is probably that in the second
configuration the rod moves back and forth in the center of the device, so the
taffy would sometimes be unstretched.  In the first configuration the rod
resolutely traverses the center of the device in a single direction each time,
leading to uniform stretching.  This is related to motions that remain `pulled
tight' as the rods move \autocite{Tumasz2013,Tumasz2013b}.  As far as the less
thorough mixing he mentions is concerned, in one turn of the handle the first
configuration gives a dilatation of~$4.2361$, while the second has~$2.4229$.
However, the second design has a larger dilatation for a full period of the
rod motion, as given in Table~\ref{tab:taffypull}.  This illustrates the
difficulties involved in comparing the efficiency of different devices.  In
its first configuration the device has a quadratic dilatation, the largest
root of~$x^2-18x+1$.  In its second configuration the dilatation is a quartic
number, the largest root of~$x^4-36x^3+5x^2-36x+1$.

\subsection*{A peculiar dilatation.}
\begin{figure}
\begin{center}
\subfigure[]{
  \includegraphics[width=.5\textwidth]{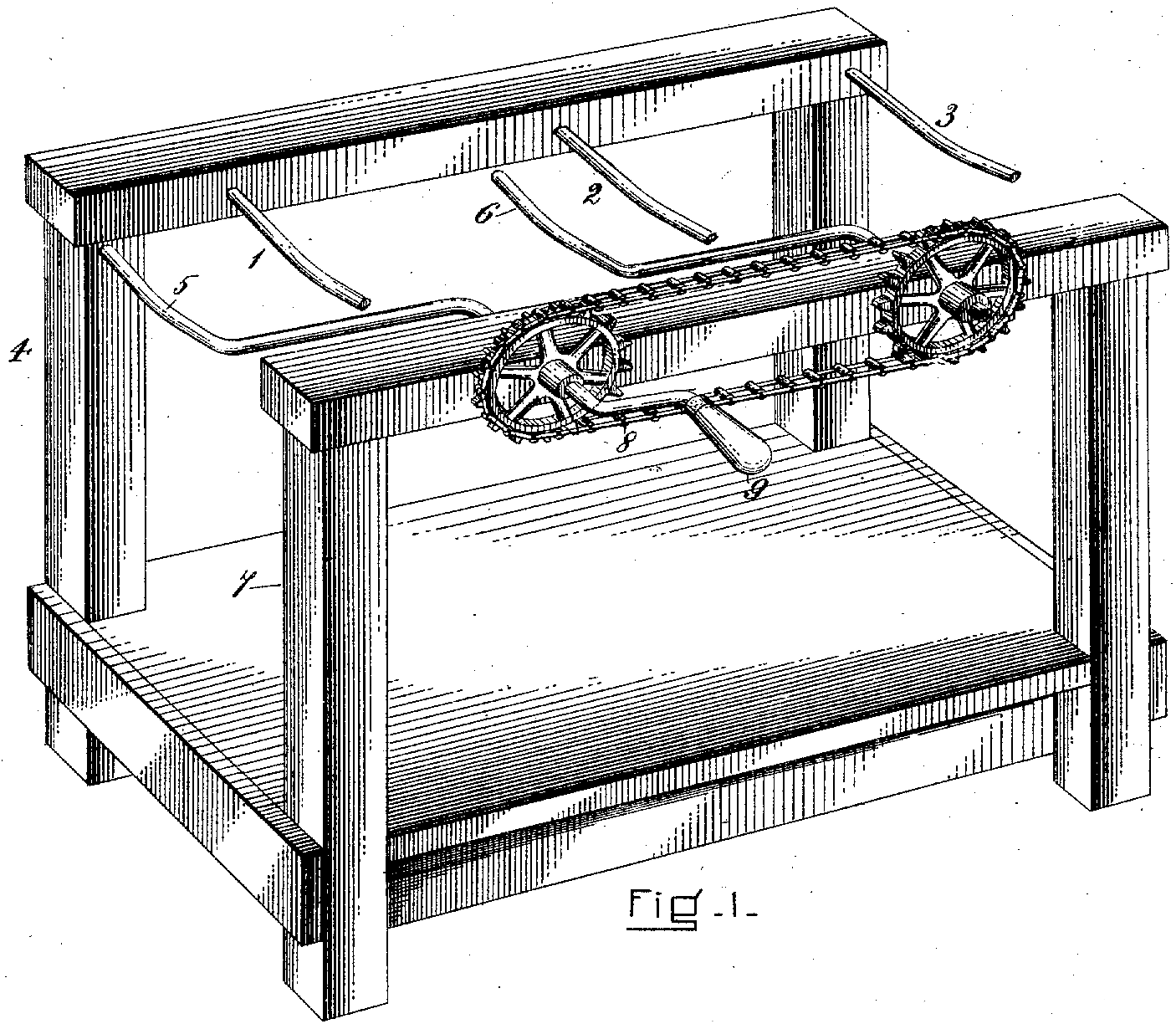}
  \label{fig:Jenner1905_device}
}\hspace{1.7em}
\subfigure[]{
  \raisebox{3em}{%
  \includegraphics[width=.4\textwidth]{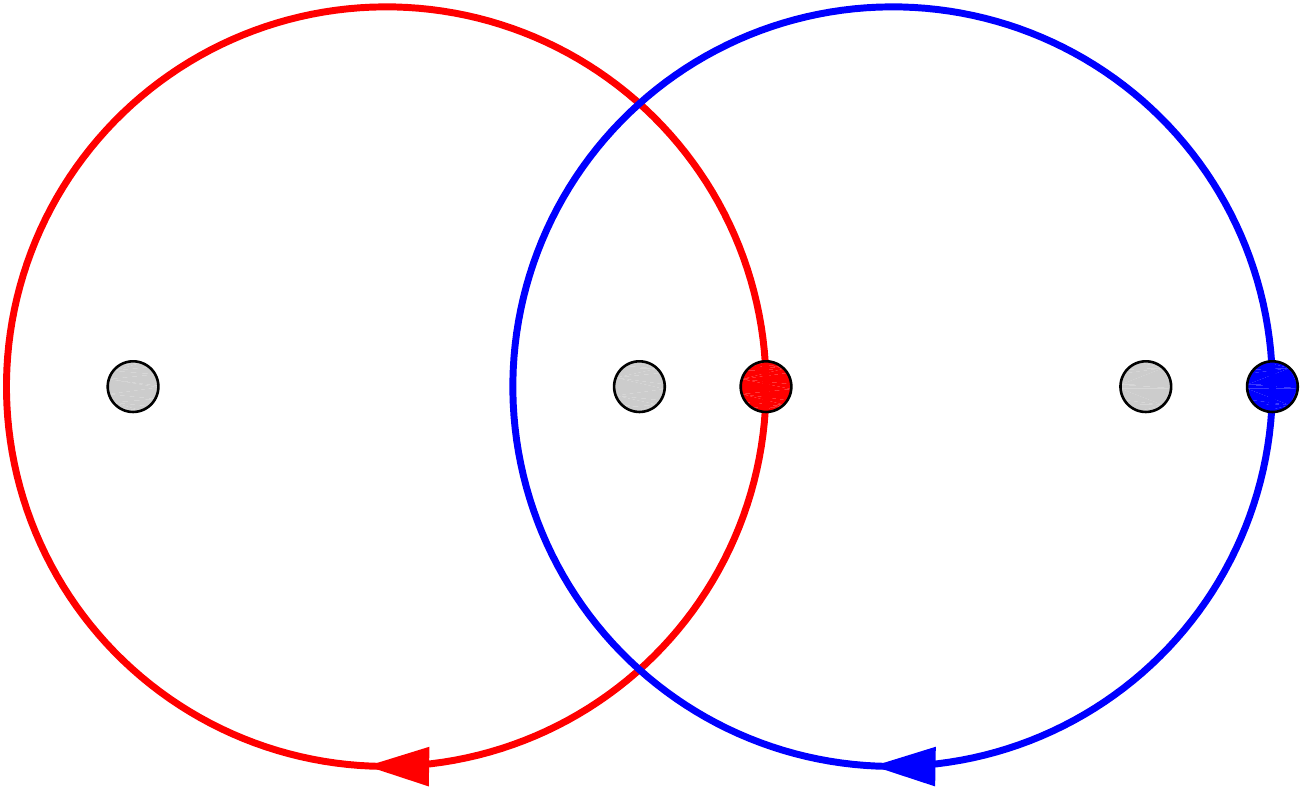}}
  \label{fig:Jenner1905_rodmotion}
}
\end{center}
\caption{Taffy puller from the patent of \textcite{jenner_candy-pulling_1905}.
  (b) The motion of the rods, with three fixed rods in gray.}
\label{fig:Jenner1905}
\end{figure}

The design of \textcite{jenner_candy-pulling_1905}, shown in
Fig.~\ref{fig:Jenner1905}, is a fairly straightforward variant of the other
devices we've seen.  From our point of view it has a peculiar property: its
dilatation is the largest root of the polynomial~$x^4-8x^3-2x^2-8x+1$, which
is the strange number~$(\varphi+\sqrt{\varphi})^2$, where~$\varphi$ is the
Golden Ratio.

\subsection*{Interlocking combs.}
\begin{figure}
\begin{center}
\subfigure[]{
  \includegraphics[width=.6\textwidth]{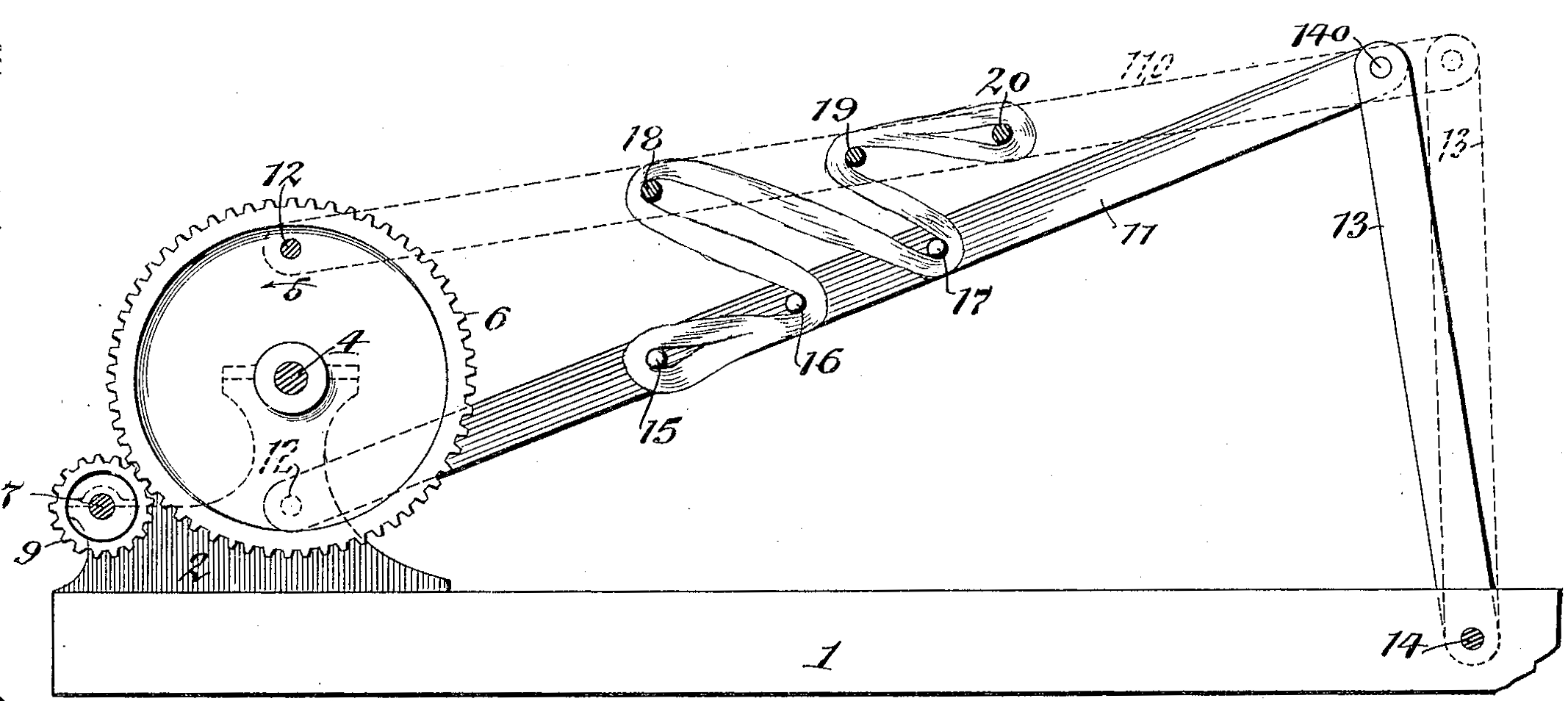}
  \label{fig:Shean1914_device}
}\hspace{.5em}
\subfigure[]{
  \raisebox{2em}{
  \includegraphics[width=.34\textwidth]{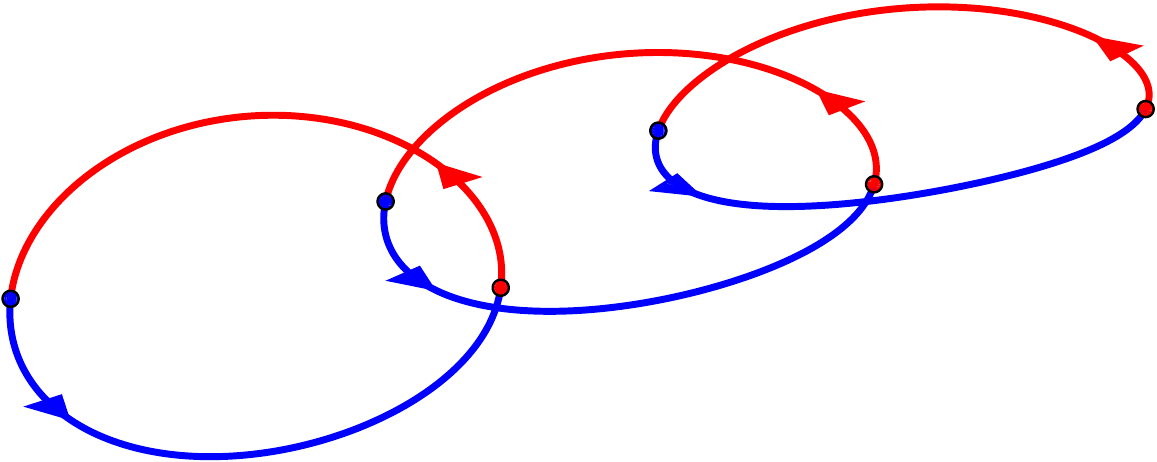}}
  \label{fig:Shean1914_rodmotion}
}
\end{center}
\caption{(a) Taffy puller from the patent of
  \textcite{shean_candy-pulling_1914}.  (b) The rod motion.}
\label{fig:Shean1914}
\end{figure}

The taffy puller of \textcite{shean_candy-pulling_1914} is shown in
Fig.~\ref{fig:Shean1914}.  The design is somewhat novel, since it is not based
directly on gears.  It consists of two interlocking `combs' of three rods
each, for a total of six moving rods.  Mathematically, this device has exactly
the same dilatation as the earlier 6-rod design (Fig.~\ref{fig:6rods}).
\keepnote{Dilatation converges to what constant?  Do periodic case.  But it's
  not quite a periodic design, though it could easily be made to be.}  A
similar comb design was later used in a device for homogenizing molten glass
\autocite{russell_apparatus_1951}.

\subsection*{A baroque design.}
\begin{figure}
\begin{center}
\subfigure[]{
  \includegraphics[height=.35\textheight]{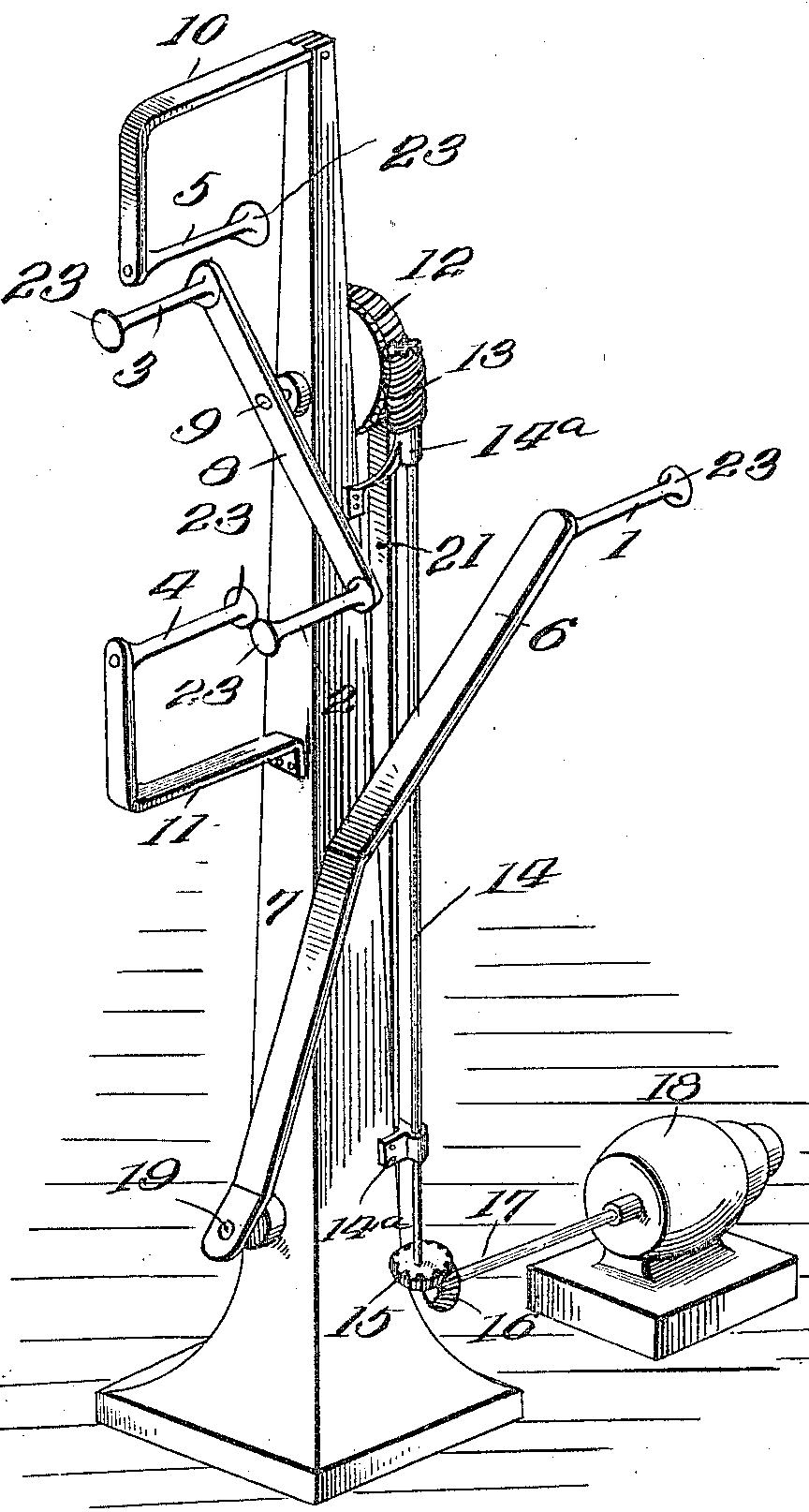}
  \label{fig:McCarthy1915_device}
}\hspace{1em}
\subfigure[]{
  \raisebox{.08\textheight}{
    \includegraphics[angle=90,height=.2\textheight]{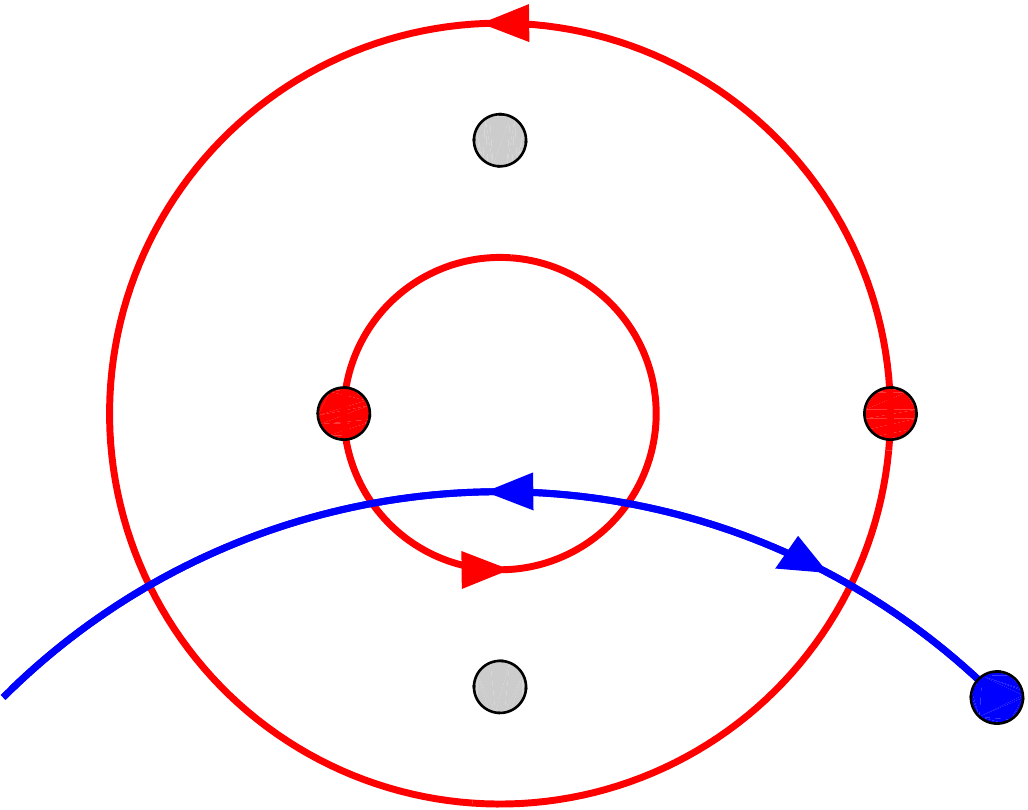}
  }
  \label{fig:McCarthy1915_rodmotion}
}
\end{center}
\caption{(a) Taffy puller from the patent of
  \textcite{mccarthy_candy-pulling_1915}.  (b) Rod motion.}
\label{fig:McCarthy1915}
\end{figure}

We finish with the intriguing design of
\textcite{mccarthy_candy-pulling_1915}, shown in Fig.~\ref{fig:McCarthy1915}.
This is the most baroque design we've encountered: it contains an oscillating
arm, rotating rods, and fixed rods.  The inventors did seem to know what they
were doing with this complexity: its dilatation is enormous at approximately
$21.2667$, the largest root of $x^4 - 20x^3 - 26x^2 - 20x + 1$.

\subsection*{Why so many designs?}
There are actually quite a few more patents for taffy pullers that were not
shown here (only U.S. patents were searched).  An obvious question is: why so
many?  Often the answer is that a new patent is created to get around an
earlier one, but the very first patents had lapsed by the 1920s and yet more
designs were introduced, so this is only a partial answer.  Perhaps there is a
natural response when looking at a taffy puller to think that we can design a
better one, since the basic idea is so simple.  At least mathematics provides
a way of making sure that we've thoroughly explored all designs, and to gauge
the effectiveness of existing ones.

%% END_ARXIV

\end{document}